\documentclass[11pt]{article}

 \usepackage{amssymb}

 \usepackage[dvips]{color}

 \newcommand{\qed}{\hfill\rule{2mm}{3mm}\vspace{4mm}}

 \newtheorem{lemma}{Lemma}[section]
 \newtheorem{theorem}{Theorem}[section]
 \newtheorem{proposition}{Proposition}[section]
 \newtheorem{corollary}{Corollary}[section]

 \def\blemma{\begin{lemma}}\def\elemma{\end{lemma}}
 \def\btheorem{\begin{theorem}}\def\etheorem{\end{theorem}}
 \def\bproposition{\begin{proposition}}\def\eproposition{\end{proposition}}
 \def\bcorollary{\begin{corollary}}\def\ecorollary{\end{corollary}}

 \def\mbb{\mathbb}\def\mbf{\mathbf}

 \def\<{\langle}\def\>{\rangle}

 \def\proof{\noindent{\it Proof.~~}}

\begin{document}

\noindent{Published in: \textit{Acta Applicandae Mathematicae.}
\textbf{88} (2005), 2: 143-175}

\bigskip\bigskip

\noindent{\LARGE\bf Conditional Log-Laplace Functionals of}

\medskip

\noindent{\LARGE\bf Immigration Superprocesses with Dependent}

\medskip

\noindent{\LARGE\bf Spatial Motion}

\bigskip

\noindent{Zenghu Li\,$^a$\footnote{ Supported by the NSFC
(No.\,10121101 and No.\,10131040).}, Hao Wang\,$^b$\footnote{
Supported by the research grant of UO.} and Jie
Xiong\,$^c$\footnote{ Research supported partially by NSA and by
Alexander von Humboldt Foundation.}}

\smallskip

\noindent{\small $^a$ School of Mathematical Sciences, Beijing
Normal University, Beijing 100875, P.R. China}

\noindent{\small E-mail: \tt lizh@email.bnu.edu.cn}

\smallskip

\noindent{\small $^b$ Department of Mathematics, University of
Oregon, Eugene OR 97403-1222, U.S.A.}

\noindent{\small E-mail: \tt haowang@darkwing.uoregon.edu}

\smallskip

\noindent{\small $^c$ Department of Mathematics, University of
Tennessee, Knoxville, TN 37996-1300, U.S.A. and}

\noindent{\small Department of Mathematics, Hebei Normal
University, Shijiazhuang 050016, P.R. China}

\noindent{\small E-mail: \tt jxiong@math.utk.edu}

\bigskip

\noindent{(Received: 7 August 2003; in finial form: 23 February
2005)}

\bigskip

\noindent{\bf Abstract.} A non-critical branching immigration
superprocess with dependent spatial motion is constructed and
characterized as the solution of a stochastic equation driven by a
time-space white noise and an orthogonal martingale measure. A
representation of its conditional log-Laplace functionals is
established, which gives the uniqueness of the solution and hence
its Markov property. Some properties of the superprocess including
an ergodic theorem are also obtained.

\bigskip


\bigskip


\bigskip\bigskip


\section{Introduction}

\setcounter{equation}{0}

A class of superprocesses with dependent spatial motion (SDSM)
over the real line $\mbb{R}$ were introduced and constructed in
Wang \cite{W97, W98}. A generalization of the model was then given
in Dawson \textit{et al} \cite{DLW01}. Let $c \in C^2_b(\mbb{R})$
and $h \in C^2_b(\mbb{R})$ and assume both $h$ and $h^\prime$ are
square-integrable. Let
 \begin{eqnarray*}
\rho(x) = \int_{\mbb{R}}h(y-x)h(y) dy, \quad x\in\mbb{R},
 \end{eqnarray*}
and $a(x) = c(x)^2 + \rho(0)$. Let $\sigma\in C^2_b(\mbb{R})^+$ be
a strictly positive function. We denote by $M(\mbb{R})$ the space
of finite Borel measures on $\mbb{R}$ endowed with a metric
compatible with its topology of weak convergence. For $f\in
C_b(\mbb{R})$ and $\mu\in M(\mbb{R})$ set $\<f,\mu\> = \int
fd\mu$. Then an SDSM $\{X_t: t\ge0\}$ is characterized by the
following martingale problem: For each $\phi \in C^2_b(\mbb{R})$,
 \begin{eqnarray}\label{1.1}
M_t(\phi) = \<\phi,X_t\> - \<\phi,X_0\>
- \frac{1}{\,2\,} \int_0^t \<a\phi^{\prime\prime},X_s\> ds,
\quad t\ge0,
 \end{eqnarray}
is a continuous martingale with quadratic variation process
 \begin{eqnarray}\label{1.2}
\<M(\phi)\>_t = \int_0^t\<\sigma\phi^2,X_s\> ds + \int_0^t
ds\int_{\mbb{R}} \<h(z - \cdot) \phi^\prime, X_s\>^2 dz.
 \end{eqnarray}
Clearly, the SDSM reduces to a usual critical branching
Dawson-Watanabe superprocess if $h(\cdot) \equiv 0$; see e.g.\
Dawson \cite{D93}. A general SDSM arises as the weak limit of
critical branching particle systems with dependent spatial motion.
Consider a family of independent Brownian motions $\{B_i(t):
t\ge0, i=1,2,\cdots\}$, the individual noises, and a time-space
white noise $\{W_t(B): t\ge 0, B\in {\cal B}(\mbb{R})\}$, the
common noise. The migration of a particle in the approximating
system with label $i$ is defined by the stochastic equation
 \begin{eqnarray}\label{1.3}
dx_i(t) = c(x_i(t)) dB_i(t) + \int_{\mbb{R}} h(y-x_i(t))W(dt,dy),
 \end{eqnarray}
where $W(ds,dy)$ denotes the time-space stochastic integral
relative to $\{W_t(B)\}$. The SDSM possesses properties very
different from those of the usual Dawson-Watanabe superprocess.
For example, a Dawson-Watanabe superprocess in $M(\mbb{R})$ is
usually absolutely continuous whereas the SDSM with $c(\cdot)
\equiv 0$ is purely atomic; see \cite{KS88} and \cite{W97, W02},
respectively.

In this paper, we consider a further extension of the model of
Wang \cite{W97, W98}. Let $b\in C^2_b(\mbb{R})$ and let $m\in
M(\mbb{R})$. A modification of the above martingale problem is to
replace (\ref{1.1}) by
 \begin{eqnarray}\label{1.4}
M_t(\phi) = \<\phi,X_t\> - \<\phi,X_0\> - t\<\phi,m\>
- \frac{1}{2}\int_0^t\<a\phi^{\prime\prime}, X_s\> ds
+ \int_0^t\<b\phi,X_s\> ds.
 \end{eqnarray}
We shall prove that there is indeed a solution $\{X_t: t\ge0\}$ to
the martingale problem given by (\ref{1.2}) and (\ref{1.4}). The
process $\{X_t: t\ge0\}$ may be regarded as a non-critical
branching \textit{SDSM with immigration} (SDSMI), where $b(\cdot)$
is the linear growth rate and $m(dx)$ gives the immigration rate.
This modification is related to the recent work of Dawson and Li
\cite{DL03}, where an interactive immigration given by
 \begin{eqnarray}\label{1.5}
\int_0^t\<q(\cdot,X_s)\phi, m\>ds
 \end{eqnarray}
was considered, where $q(\cdot,\cdot)$ is a function on
$\mbb{R}\times M(\mbb{R})$ representing a state dependent
immigration density. However, it was assumed in \cite{DL03} that
$b(\cdot) \equiv c(\cdot) \equiv 0$ and the approach there relies
essentially on the purely atomic property of the process, which is
not available for the present model.

The main purpose of the paper is to give a representation of the
conditional log-Laplace functionals of solution of (\ref{1.2}) and
(\ref{1.4}) and to illustrate some applications of the
representation. This approach was stimulated by Xiong \cite{X04},
who established a similar characterization for the model of
Skoulakis and Adler \cite{SA01}. The key idea of the
representation is to decompose the martingale (\ref{1.4}) into two
orthogonal components, which arise respectively from the migration
and the branching. Since the decomposition uses additional
information which is not provided by (\ref{1.2}) and (\ref{1.4}),
we shall start with the corresponding particle system and consider
the high density limit following \cite{DVW00}. In this way, we can
easily separate the two kinds of noises. It turns out that the
common migration noise $\{W(ds,dy)\}$ remains after the limit
procedure and the limit process satisfies the following martingale
problem: For each $\phi \in C^2_b(\mbb{R})$,
 \begin{eqnarray}\label{1.6}
Z_t(\phi)
&=& \<\phi,X_t\> - \<\phi,X_0\> - t\<\phi,m\>
- \frac{1}{2}\int_0^t\<a\phi^{\prime\prime}, X_s\> ds  \nonumber  \\
& & + \int_0^t\<b\phi,X_s\> ds -
\int_0^t\int_{\mbb{R}}\<h(y-\cdot) \phi^\prime,X_s\> W(ds,dy)
 \end{eqnarray}
is a continuous martingale orthogonal to $\{W_t(\phi)\}$ with
quadratic variation process
 \begin{eqnarray}\label{1.7}
\<Z(\phi)\>_t
=
\int_0^t\<\sigma\phi^2, X_s\>ds.
 \end{eqnarray}
This formulation suggests that we may regard $\{X_t: t\ge0\}$
as a generalized inhomogeneous Dawson-Watanabe superprocess with
immigration, where
 \begin{eqnarray*}
\int_{\mbb{R}}h(y-\cdot) W(dt,dy)
 \end{eqnarray*}
gives a generalized drift in the underlying migration. Based on
the techniques developed in Kurtz and Xiong \cite{KX99, X04}, we
prove that for each $\phi \in H_1(\mbb{R})\cap C_b(\mbb{R})$ there
is a pathwise unique solution of the non-linear SPDE
 \begin{eqnarray}\label{1.8}
\psi_{r,t}(x)
&=&
\phi(x) + \int_r^t \bigg[\frac{1}{2} a(x)\psi_{s,t}
^{\prime\prime}(x) - \frac{1}{2}\sigma(x)
\psi_{s,t}(x)^2\bigg] ds \nonumber \\
& & - \int_r^t b(x)\psi_{s,t}(x) ds + \int_r^t\int_{\mbb{R}}
h(y-x)\psi_{s,t}^\prime(x) \cdot W(ds,dy),
 \end{eqnarray}
where the last term on the right hand side denotes the backward
stochastic integral with respect to the white noise. Then we show
that the conditional log-Laplace functionals of $\{X_t: t\ge0\}$
given $\{W(ds,dy)\}$ can be represented by the solution of
(\ref{1.8}). The representation of the conditional log-Laplace
functionals is proved by direct analysis based on (\ref{1.6}),
(\ref{1.7}) and (\ref{1.8}). This approach is different from that
of Xiong \cite{X04}, where a Wong-Zakai type approximation was
used. The idea of conditional log-Laplace approach has also been
used by Crisan \cite{C04} for a different model. In fact, the
approach in Section 5 is adapted from \cite{C04} which simplifies
our original arguments. It is well-known that non-conditional
log-Laplace functionals play very important roles in the study of
classical Dawson-Watanabe superprocesses.

We shall see that conditional Laplace functionals are almost as
efficient as the non-conditional Laplace functionals in studying
some properties of the SDSMI. In particular, the characterization
of the conditional Laplace functionals gives immediately the
uniqueness of solution of (\ref{1.6}) and (\ref{1.7}), which in
turn implies the Markov property of $\{X_t: t\ge0\}$. It follows
that $\{X_t: t\ge0\}$ is a diffusion process with generator ${\cal
L}$ given by
 \begin{eqnarray}\label{1.9}
{\cal L} F(\mu)
 &=&
\frac{1}{2}\int_{\mbb{R}^2}\rho(x-y)\frac{d^2}{dxdy}
\frac{\delta^2
F(\mu)}{\delta\mu(x)\delta\mu(y)} \mu(dx)\mu(dy)   \nonumber \\
 & &
+ \frac{1}{2}\int_{\mbb{R}} a(x)\frac{d^2}{dx^2} \frac{\delta
F(\mu)}{\delta\mu(x)}\mu(dx) + \frac{1}{2} \int_{\mbb{R}}
\sigma(x)\frac{\delta^2 F(\mu)} {\delta\mu(x)^2} \mu(dx) \nonumber \\
 & &
- \int_{\mbb{R}} b(x)\frac{\delta F(\mu)} {\delta\mu(x)} \mu(dx) +
\int_{\mbb{R}} \frac{\delta F(\mu)} {\delta\mu(x)} m(dx),
 \end{eqnarray}
where
 \begin{eqnarray}\label{1.10}
\frac{\delta F(\mu)}{\delta\mu(x)} = \lim_{r\to
0^+}\frac{1}{\,r\,}[F(\mu + r\delta_x) - F(\mu)]
 \end{eqnarray}
and $\delta^2F(\mu) / \delta\mu(x)\delta\mu(y)$ is defined in the
same way with $F$ replaced by $(\delta F/ \delta\mu(y))$ on the
right hand side; see Section 3. We also prove some properties of
the SDSMI including an ergodic theorem. There are also some other
applications of the conditional log-Laplace functional. For
instance, based on this characterization the conditional excursion
theory of the SDSM have been developed in \cite{LWX04b}. However,
consideration of the interactive immigration (\ref{1.5}) for this
present process seems sophisticated.

The remainder of the paper is organized as follows. In Section 2
we give a formulation of the system of branching particles with
dependent spatial motions and immigration. Some useful estimates
of the moments of the system are also given. In Section 3 we
obtain a solution of the martingale problem (\ref{1.6}) and
(\ref{1.7}) as the high density limit of a sequence of particle
systems. The existence and uniqueness of the solution of
(\ref{1.8}) is established in Section 4. In Section 5 we give the
representation of the conditional log-Laplace functionals of the
solution of (\ref{1.6}) and (\ref{1.7}). Some properties of the
SDSMI are discussed in Section 6.


\section{Branching particle systems}

\setcounter{equation}{0}

The main purpose of this section is to give an explicit
construction for the immigration branching particle system with
dependent spatial motion by modifying the constructions of
\cite{DVW00, W86}. This construction provides a useful set up of
the process.

We start with a simple interacting particle system. Let $\theta >
0$ be a constant and $(c,h)$ be given as in the introduction. Let
$N(\mbb{R}) \subset M(\mbb{R})$ be the set of integer-valued
measures on $\mbb{R}$ and let $M_\theta(\mbb{R}) := \{\theta^{-1}
\sigma: \sigma \in N(\mbb{R})\}$. Given $\{a_i: i=1,\cdots,n\}$,
let $\{x_i(t): t\ge0, i=1,\cdots,n\}$ be given by
 \begin{eqnarray}\label{2.1}
x_i(t) = a_i + \int_0^tc(x_i(s)) dB_i(s) + \int_0^t\int_{\mbb{R}}
h(y-x_i(s))W(dy, ds).
 \end{eqnarray}
We may define a measure-valued process $\{X_t: t\ge0\}$ by
 \begin{eqnarray}\label{2.2}
 \<\phi,X_t\>  = \sum_{i=1}^n \theta^{-1}\phi(x_i(t)), \qquad t\ge0.
 \end{eqnarray}
By the discussions in \cite{DLW01, W97, W98}, the process $\{X_t:
t\ge0\}$ is a diffusion in $M_\theta (\mbb{R})$. Let ${\cal
A}_\theta$ denote the generator of this diffusion process. If
$F_{f,\{\phi_i\}} (\mu) := f(\<\phi_1,\mu\>, \cdots,
\<\phi_n,\mu\>)$ for $f\in C^2_0(\mbb{R}^n)$ and
$\{\phi_i\}\subset C^2_b(\mbb{R})$, by It\^o's formula it is easy
to see that
 \begin{eqnarray}\label{2.3}
{\cal A}_\theta F_{f,\{\phi_i\}}(\mu) &=&
\frac{1}{2}\sum_{i,j=1}^n f_{ij}^{\prime\prime}(\<\phi_1,
\mu\>,\cdots,\<\phi_n,\mu\>)\int_{\mbb{R}^2}\rho(x-y)
\phi_i^\prime(x)\phi_j^\prime(y)\mu(dx)\mu(dy)  \nonumber \\
& & +\,\frac{1}{2}\sum_{i=1}^n f_i^\prime(\<\phi_1,\mu\>,\cdots,
\<\phi_n,\mu\>)\<a\phi_i^{\prime\prime},\mu\> \nonumber \\
& & +\,\frac{1}{2\theta}\sum_{i,j=1}^n f_{ij}^{\prime\prime}
(\<\phi_1,\mu\>,\cdots,\<\phi_n,\mu\>) \<c^2\phi_i^\prime
\phi_j^\prime,\mu\>.
 \end{eqnarray}
More generally, if $F$ is a function on $M_\theta(\mbb{R})$ that
can be extended to a sufficiently smooth function on $M(\mbb{R})$,
then
 \begin{eqnarray}\label{2.4}
{\cal A}_\theta F(\mu) &=&
\frac{1}{2}\int_{\mbb{R}^2}\rho(x-y)\frac{d^2}{dxdy}
\frac{\delta^2 F(\mu)}{\delta\mu(x)\delta\mu(y)}
\mu(dx)\mu(dy) \nonumber  \\
& & +\,\frac{1}{2}\int_{\mbb{R}} a(x)\frac{d^2}{dx^2}\frac{\delta
F(\mu)}{\delta\mu(x)}\mu(dx) \nonumber  \\
& & +\,\frac{1}{2\theta}\int_{\mbb{R}^2} c(x)c(y) \frac{d^2}{dxdy}
\frac{\delta^2 F(\mu)}{\delta\mu(x)\delta\mu(y)}
\delta_x(dy)\mu(dx),
 \end{eqnarray}
where $\delta F(\mu) / \delta\mu(x)$ and $\delta^2F(\mu) /
\delta\mu(x)\delta\mu(y)$ are defined as in the introduction. This
can be seen by approximating the function $F$ by functions of the
form $F_{f,\{\phi_i\}}$.

A more interesting particle system involves branching and
immigration. Let $\gamma>0$ be a constant and let $m\in
M(\mbb{R})$. Let $p(x,\cdot) = \{p_0(x), p_1(x), p_2(x), \cdots\}$
be a family of discrete probability distributions which measurably
depends on the index $x\in \mbb{R}$ and satisfies $p_1(\cdot)
\equiv 0$. In addition, we assume that
 \begin{eqnarray}\label{2.5}
q(x) := \sum_{i=1}^\infty ip_i(x), \quad x\in \mbb{R},
 \end{eqnarray}
is a bounded function. We shall construct an immigration branching
particle system with parameters $(a,\rho,\gamma,p,\theta
m,1/\theta)$.

Let ${\cal A}$ be the set of all strings of the form $\alpha =
n_0n_1 \cdots n_{l(\alpha)}$, where $l(\alpha)$ is the length of
$\alpha$ and the $n_j$ are non-negative integers with $0\le n_0\le
1$ and $n_j\ge1$ for $j\ge1$. We shall label the particles by the
strings in ${\cal A}$. We here use the first digit $n_0$ in the
string to distinguish the aboriginal and the immigratory
particles. More precisely, strings started with $0$ refer to
descendants of aboriginal ancestors and strings started with $1$
refer to descendants of immigratory ancestors. (Note that the
first digit is not counted in the length $l(\alpha)$.) We provide
${\cal A}$ with the arboreal ordering, that is, $m_0\cdots m_p
\prec n_0\cdots n_q$ if and only if $p\le q$ and $m_0=n_0, \cdots,
m_p=n_p$. Then $\alpha$ has exactly $l(\alpha)$ predecessors,
which we denote respectively by $\alpha-1$, $\alpha-2$, $\cdots$,
$\alpha - l(\alpha)$. For example, if $\alpha = 12431$, then
$\alpha-2 = 124$ and $\alpha-4 = 1$.

We need a collection of random variables to construct the
immigration branching particle system. Let $\{a_{01}, \cdots,
a_{0n}\}$ be a finite sequence of real-valued random variables.
Let $\{W(ds,dx): s\ge0, x\in \mbb{R}\}$ be a time-space white
noise and $\{N(ds,dx): s\ge0, x\in \mbb{R}\}$ a Poisson random
measure with intensity $\theta dsm(dx)$. We shall assume $\<1,m\>
>0$, otherwise the construction of the immigration part is
trivial. In this case, we can enumerate the atoms of $N(ds,dx)$ as
 \begin{eqnarray}\label{2.6}
\{(s_i,a_{1i}): 0<s_1<s_2<\cdots, a_{1i}\in\mbb{R}\}.
 \end{eqnarray}
We also define the families
 \begin{eqnarray}\label{2.7}
\{B_\alpha(t): t\ge0, \alpha\in{\cal A}\},\quad \{S_\alpha:
\alpha\in{\cal A}\},\quad \{\eta_{a,\alpha}:
a\in\mbb{R},\alpha\in{\cal A}\},
 \end{eqnarray}
where $\{B_\alpha\}$ are independent standard Brownian motions,
$\{S_\alpha\}$ are i.i.d.\ exponential random variables with
parameter $\gamma$, and $\{\eta_{a,\alpha}\}$ are independent
random variables with distribution $p(a,\cdot)$. We assume that
the families $\{W(ds,dx)\}$, $\{N(ds,dx)\}$, $\{a_{0i}\}$,
$\{B_\alpha\}$, $\{S_\alpha\}$ and $\{\eta_{a,\alpha}\}$ are
independent.

We define $\beta_{0n_1} = 0$ if $1\le n_1 \le n$ and $\beta_{0n_1}
= \infty$ if $n_1 > n$, and define $\beta_{1n_1} = s_{n_1}$ for
all $n_1\ge 1$. For $\alpha\in{\cal A}$ with $l(\alpha)= 1$ we let
$\zeta_{\alpha} = \beta_{\alpha} + S_{\alpha}$. Heuristically,
$S_\alpha$ is the life-span of the particle with label $\alpha$,
$\beta_\alpha$ is its birth time and $\zeta_\alpha$ is its death
time. The random variables $a_\alpha$ defined above can be
interpreted as the birth place of the particle with label
$\alpha$. The trajectory $\{x_\alpha(t): t\ge\beta_\alpha\}$ of
the particle is the solution of the equation
 \begin{eqnarray}\label{2.8}
x(\beta_\alpha+t) = a_\alpha +
\int_{\beta_\alpha}^{\beta_\alpha+t}c(x(s))dB_\alpha(s) +
\int_{\beta_\alpha}^{\beta_\alpha+t}\int_{\mbb{R}} h(y-x(s))
W(ds,dy).
 \end{eqnarray}
For $\alpha\in {\cal A}$ with $l(\alpha)>1$ the trajectory
$\{x_\alpha (t): t\ge\beta_\alpha\}$ is defined by the above
equation with $a_\alpha = x_{\alpha-1} (\zeta_{\alpha-1}^-)$,
$\zeta_{\alpha} = \beta_{\alpha} + S_{\alpha}$ and
 \begin{eqnarray}\label{2.9}
\beta_\alpha
=\left\{\begin{array}{ll}
\zeta_{\alpha-1}   &\mbox{ if $n_{l(\alpha)} \le
\eta_{x_{\alpha-1}(\zeta_{\alpha-1}-),\alpha-1}$}  \\
\infty             &\mbox{ if $n_{l(\alpha)} >
\eta_{x_{\alpha-1}(\zeta_{\alpha-1}-),\alpha-1}$,}
\end{array}\right.
 \end{eqnarray}
where $x_{\alpha-1}(\zeta_{\alpha-1}-)$ denotes the left limit of
$x_{\alpha-1}(t)$ at $t=\zeta_{\alpha-1}$. Clearly,
 \begin{eqnarray}\label{2.10}
\<\phi,Y_t\> = \sum_{\alpha\in{\cal A}} \theta^{-1}
\phi(x_\alpha(t)) 1_{[\beta_\alpha,\zeta_\alpha)}(t), \qquad
t\ge0.
 \end{eqnarray}
defines an $M_\theta(\mbb{R})$-valued process $\{Y_t: t\ge0\}$. It
is easy to see that $\{Y_t: t\ge0\}$ has countably many jumps, and
between those jumps it behaves just as the diffusion process
$\{X_t: t\ge0\}$ constructed by (\ref{2.2}). We call $\{Y_t:
t\ge0\}$ an \textit{immigration branching particle system} with
parameters $(c,h,\gamma,p,\theta m,1/\theta)$. Intuitively,
$p(x,\cdot)$ gives the location dependent offspring distribution
and $\{N(ds,dx)\}$ gives the landing times and sites of the
immigrants.

Indeed, we may regard $\{Y_t: t\ge0\}$ as a concatenation of a
sequence of independent copies of $\{X_t: t\ge0\}$. We refer the
reader to \cite{S88} for discussions of concatenation of general
Markov processes. As in \cite{LLW04} it can be seen that $\{Y_t:
t\ge0\}$ is a Markov process with generator ${\cal L}_\theta :=
{\cal A}_\theta + {\cal B}_\theta$, where
 \begin{eqnarray}\label{2.11}
{\cal B}_\theta F(\mu) &=& \sum_{j=0}^\infty \int_{\mbb{R}} \theta
\gamma p_j(x) \big[F\big(\mu + (j-1)\theta^{-1} \delta_x\big)
- F(\mu)\big] \mu(dx)  \nonumber \\
& &\quad + \int_{\mbb{R}} \theta \big[F\big(\mu + \theta^{-1}
\delta_x\big) - F(\mu)\big] m (dx).
 \end{eqnarray}
The first term on the right hand side of (\ref{2.11}) represents
the jumps given by the branching and the second terms represents
the jumps given by the immigration. In particular, it is easy to
show that
 \begin{eqnarray}\label{2.12}
{\cal B}_\theta F_{f,\{\phi_i\}}(\mu) &=& \sum_{j=0}^\infty
\int_{\mbb{R}} \theta \gamma p_j(x)
\big[f(\<\phi_1,\mu\>+\theta^{-1}\phi_1(x),\cdots,
\<\phi_n,\mu\>+\theta^{-1}\phi_n(x))  \nonumber \\
& &\hskip2cm - f(\<\phi_1,\mu\>,\cdots,
\<\phi_n,\mu\>)\big] \mu(dx)  \nonumber \\
& & + \int_{\mbb{R}} \theta
\big[f(\<\phi_1,\mu\>+\theta^{-1}\phi_1(x),\cdots,
\<\phi_n,\mu\>+\theta^{-1}\phi_n(x))  \nonumber \\
& &\hskip2cm - f(\<\phi_1,\mu\>,\cdots, \<\phi_n,\mu\>)\big] m
(dx).
 \end{eqnarray}
Let ${\cal D}_1 ({\cal L}_\theta)$ denote the collection of all
functions $F_{f,\{\phi_i\}}$ with $f\in C^2_0(\mbb{R}^n)$ and
$\{\phi_i\}\subset C^2_b(\mbb{R})$. By the general theory of
Markov processes, we have the following

\btheorem\label{t2.1} The process $\{Y_t: t\ge0\}$ defined by
(\ref{2.10}) solves the $({\cal L}_\theta,{\cal D}_1 ({\cal
L}_\theta))$-martingale problem, that is, for each $F \in {\cal
D}_1 ({\cal L}_\theta)$,
 \begin{eqnarray*}
F(X_t) - F(X_0) - \int_0^t {\cal L}_\theta F(X_s)ds, \qquad t\ge0,
 \end{eqnarray*}
is a martingale. \etheorem

Let us give another useful formulation of the immigration particle
system. {From} (\ref{2.8}), (\ref{2.10}) and It\^o's formula we
get
 \begin{eqnarray*}
\<\phi,Y_t\> &=& \<\phi,Y_0\>
+ \sum_{i=1}^\infty \theta^{-1}\phi(a_{1i})1_{(0,t]}(s_i)  \\
& & + \sum_{\alpha\in{\cal A}} [\eta_{x_\alpha(\zeta_\alpha-),
\alpha} - 1]\theta^{-1}\phi(x_\alpha(\zeta_\alpha-))
1_{(0,t]}(\zeta_\alpha)  \\
& & + \sum_{\alpha\in{\cal A}}\int_0^t\theta^{-1}\phi^\prime
(x_\alpha(s))1_{[\beta_\alpha,\zeta_\alpha)}(s)c(x_\alpha(s))
dB_\alpha(s)  \\
& & + \sum_{\alpha\in{\cal A}}\int_0^t\int_{\mbb{R}}\theta^{-1}
\phi^\prime(x_\alpha(s))1_{[\beta_\alpha,\zeta_\alpha)}
(s) h(y-x_\alpha(s)) W(ds,dy)  \\
& & + \frac{1}{2}\sum_{\alpha\in{\cal A}}\int_0^t\theta^{-1}
\phi^{\prime\prime}(x_\alpha(s))1_{[\beta_\alpha,
\zeta_\alpha)}(s) a(x_\alpha(s))ds,
 \end{eqnarray*}
which can be rewritten as
 \begin{eqnarray}\label{2.13}
\<\phi,Y_t\> &=& \<\phi,Y_0\> + \int_{(0,t]}\int_{\mbb{R}}
\theta^{-1}\phi(x)N(ds,dx) \nonumber \\
 & &
+ \sum_{\alpha\in{\cal A}} [\eta_{x_\alpha(\zeta_\alpha-), \alpha}
- 1]\theta^{-1}\phi(x_\alpha(\zeta_\alpha-))
1_{(0,t]}(\zeta_\alpha)  \nonumber \\
 & &
+ \sum_{\alpha\in{\cal A}}\int_0^t\theta^{-1}\phi^\prime
(x_\alpha(s))1_{[\beta_\alpha,\zeta_\alpha)}(s)
c(x_\alpha(s)) dB_\alpha(s)  \nonumber  \\
 & &
+ \int_0^t\int_{\mbb{R}} \<h(y-\cdot)\phi^\prime,Y_s\> W(ds,dy) +
\frac{1}{2}\int_0^t \<a\phi^{\prime\prime}, Y_s\> ds.
 \end{eqnarray}
On the right hand side, the second term comes from the
immigration, the third term represents branching of the particles,
and the last three terms are determined by the spatial motion. It
is not hard to see that, for any $\psi\in C_b(\mbb{R})$,
 \begin{eqnarray}\label{2.14}
U_t(\psi) := \sum_{\alpha\in{\cal A}}\int_0^t \theta^{-1}
\psi(x_\alpha(s)) 1_{[\beta_\alpha,\zeta_\alpha)}(s)
c(x_\alpha(s)) dB_\alpha(s)
 \end{eqnarray}
is a continuous local martingale with quadratic variation process
 \begin{eqnarray}\label{2.15}
\<U(\psi)\>_t := \int_0^t\<\theta^{-1} c^2\psi^2,Y_s\> ds.
 \end{eqnarray}
In the sequel, we assume
 \begin{eqnarray}\label{2.16}
\sigma(x) = \sum_{i=0}^\infty p_i(x)(i-1)^2, \qquad x\in \mbb{R},
 \end{eqnarray}
is a bounded function on $\mbb{R}$.

\bproposition\label{p2.1} For any $\phi\in C_b(\mbb{R})$,
 \begin{eqnarray}\label{2.17}
Z_t(\phi) := \sum_{\alpha\in{\cal A}}
[\eta_{x_\alpha(\zeta_\alpha-),\alpha} - 1] \theta^{-1}
\phi(x_\alpha(\zeta_\alpha-))1_{(0,t]} (\zeta_\alpha) -
\int_0^t\<\gamma(q-1)\phi,Y_s\> ds
 \end{eqnarray}
is a local martingale with predictable quadratic variation process
 \begin{eqnarray}\label{2.18}
\<Z(\phi)\>_t = \int_0^t\<\theta^{-1}\gamma\sigma\phi^2,Y_s\> ds.
 \end{eqnarray}
\eproposition

\proof Recall that $\{S_\alpha\}$ are i.i.d.\ exponential random
variables with parameter $\gamma$. Let
 \begin{eqnarray}\label{2.19}
J_t(\phi) = \sum_{\alpha\in{\cal A}} \theta^{-1}[\eta_{x_\alpha
(\zeta_\alpha-),\alpha} - 1]\phi(x_\alpha(\zeta_\alpha-))
1_{(0,t]}(\zeta_\alpha).
 \end{eqnarray}
Observe that the process $\{J_t(\phi): t\ge0\}$ jumps only when a
particle in the population splits. It is not hard to show that
$\{(Y_t,J_t(\phi)): t\ge0\}$ is a Markov process with generator
${\cal J}_\theta$ such that
 \begin{eqnarray*}
{\cal J}_\theta F(\mu,z) &=& {\cal A}_\theta F(\cdot,z)(\mu) +
\int_{\mbb{R}} \theta [F(\mu + \theta^{-1}\delta_x,z)
- F(\mu,z)] m(dx) \nonumber \\
& & + \sum_{j=0}^\infty \int_{\mbb{R}} \theta \gamma p_j(x) [F(\mu
+ (j-1)\theta^{-1} \delta_x,z + (j-1)\theta^{-1}\phi(x)) -
F(\mu,z)] \mu(dx).
 \end{eqnarray*}
In particular, if $F(\mu,z) = z$, then
 \begin{eqnarray*}
{\cal J}_\theta F(\mu,z) = \sum_{j=0}^\infty \int_{\mbb{R}} \gamma
p_j(x)(j-1)\phi(x) \mu(dx) = \<\gamma(q-1)\phi,\mu\>.
 \end{eqnarray*}
This shows that (\ref{2.17}) is a local martingale. Let $\Delta_n
:= \{0= t_{n,0} < t_{n,1} < \cdots < t_{n,n} =t\}$ be a sequence
of partitions of $[0,t]$ such that $D_n := \max_{1\le i\le n}
|t_{n,i} - t_{n,i-1}| \to 0$ as $n \to \infty$. Since the second
term on the right hand side of (\ref{2.17}) is of locally finite
variations, we have
 \begin{eqnarray*}
[Z(\phi)]_{t\land\tau_l}
&:=& \lim_{n\to \infty}\sum_{i=0}^n |Z_{t_{n,i}\land\tau_l}(\phi)
- Z_{t_{i-1}\land\tau_l}(\phi)|^2 \nonumber \\
&=& \sum_{\alpha\in{\cal A}} \theta^{-2} [\eta_{x_\alpha
(\zeta_\alpha-),\alpha} - 1]^2\phi(x_\alpha (\zeta_\alpha-))^2
1_{(0,t\land\tau_l]}(\zeta_\alpha).
 \end{eqnarray*}
By martingale theory, $Z_{t\land\tau_l} (\phi)^2 - [Z(\phi)]_{t
\land \tau_l}$ is a martingale. Note that $[Z(\phi)]_{t \land
\tau_l}$ has same jump times as $J_{t\land\tau_l} (\phi)$ but with
squared jump sizes. By an argument similar to the beginning of
this proof, we conclude that $[Z(\phi)]_{t \land \tau_l} -
\<Z(\phi)\>_{t \land \tau_l}$ is a martingale. Then
$\<Z(\phi)\>_{t \land \tau_l}$ is a predictable process such that
$Z_{t\land\tau_l}(\phi)^2 - \<Z(\phi)\>_{t \land \tau_l}$ is a
martingale, implying the desired result. \qed

Let $\tilde N(ds,dx) = N(ds, dx) - \theta dsm(dx)$. Note that the
assumptions on independence imply that the four martingale
measures $\{W(ds,dx)\}$, $\{\tilde N(ds,dx)\}$, $\{Z(ds,dx)\}$ are
$\{U(ds,dx)\}$ are orthogonal to each other. Now we may rewrite
(\ref{2.13}) into
 \begin{eqnarray}\label{2.20}
\<\phi,Y_t\> &=& \<\phi,Y_0\> + t\<\phi,m\> +
\int_{(0,t]}\int_{\mbb{R}}
\theta^{-1}\phi(x) \tilde N(ds,dx)  \nonumber  \\
& & + \int_0^t\<\gamma(q-1)\phi,Y_s\> ds
+ Z_t(\phi) + U_t(\phi^\prime)  \nonumber  \\
& & + \int_0^t\int_{\mbb{R}} \<h(y-\cdot)\phi^\prime,Y_s\>
W(ds,dy) + \frac{1}{2}\int_0^t \<a\phi^{\prime\prime},Y_s\> ds.
 \end{eqnarray}
Clearly, the third term on the right hand side of (\ref{2.20}) has
a c\`adl\`ag modification. By \cite[p.69, Theorem VI.4]{DM82}, the
martingale $\{Z_t(\phi): t\ge0\}$ has a c\`adl\`ag modification.
All other terms on the right hand side have continuous
modifications. Therefore, the measure-valued process $\{Y_t:
t\ge0\}$ has a c\`adl\`ag modification and (\ref{2.20}) gives an
SPDE formulation of this immigration branching particle system.
The following result shows that (\ref{2.14}) and (\ref{2.17}) are
in fact square-integrable martingales.

\bproposition\label{p2.2} Let $B_1 := \|\gamma(q-1)\|$ and $B_2 :=
\|\theta \gamma \sigma\|$, where $\|\cdot\|$ denotes the supremum
norm. Then there is a locally bounded function $C_2$ on
$\mbb{R}_+^3$ such that
 \begin{eqnarray}\label{2.21}
\mbf{E}\{\mbox{$\sup_{0\le s\le t}$}\<1,Y_s\>^2\} \le
C_2(B_1,B_2,t) (1+\<1,\mu\>^2 + \<1,m\>^2), \quad t\ge 0.
 \end{eqnarray}
\eproposition

\proof Applying (\ref{2.20}) to $\phi \equiv 1$ we get
 \begin{eqnarray}\label{2.22}
\<1,Y_t\> = \<1,\mu\> + \theta^{-1} N((0,t] \times \mbb{R}) +
\int_0^t\<\gamma(q-1),Y_s\> ds + Z_t(1),
 \end{eqnarray}
where $N((0,t] \times \mbb{R})$ is a Poisson random variable with
parameter $\theta t\<1,m\>$ and $\{Z_t(1): t\ge0\}$ is a local
martingale with quadratic variation process
 \begin{eqnarray}\label{2.23}
\<Z(1)\>_t = \int_0^t\<\theta^{-1}\gamma\sigma,Y_s\> ds.
 \end{eqnarray}
Based on (\ref{2.22}) and (\ref{2.23}), the desired estimate
follows by an application of Gronwall's inequality. \qed


\section{Stochastic equation of the SDSMI}

\setcounter{equation}{0}

Let $(c, h, \sigma, b, m)$ be given as in the introduction.
Suppose that $W(ds,dx)$ is a time-space white noise. For $\mu\in
M(\mbb{R})$ we consider the stochastic equation:
 \begin{eqnarray}\label{3.1}
\<\phi,X_t\> &=& \<\phi,\mu\> + t\<\phi,m\>
+ \frac{1}{2}\int_0^t\<a\phi^{\prime\prime}, X_s\> ds
- \int_0^t\<b\phi,X_s\> ds \nonumber \\
& & + \int_0^t\int_{\mbb{R}}\phi(y) Z(ds,dy) +
\int_0^t\int_{\mbb{R}}\<h(y-\cdot) \phi^\prime,X_s\> W(ds,dy),
 \end{eqnarray}
where $Z(ds,dy)$ is an orthogonal martingale measure which is
orthogonal to the white noise $W(ds,dy)$ and has covariation
measure $\sigma(y)X_s(dy)ds$. Clearly, this is equivalent to the
martingale problem given by (\ref{1.6}) and (\ref{1.7}). We shall
prove that (\ref{3.1}) has a weak solution $\{X_t: t\ge0\}$, which
will serve as a candidate of the SDSMI with parameters $(c, h,
\sigma, b, m)$. For a function $F$ on $M(\mbb{R})$, let
 \begin{eqnarray}\label{3.2}
{\cal A} F(\mu) &=&
\frac{1}{2}\int_{\mbb{R}^2}\rho(x-y)\frac{d^2}{dxdy}
\frac{\delta^2 F(\mu)}{\delta\mu(x)\delta\mu(y)}
\mu(dx)\mu(dy) \nonumber  \\
& &\quad +\,\frac{1}{2}\int_{\mbb{R}} a(x)\frac{d^2}{dx^2}
\frac{\delta F(\mu)}{\delta\mu(x)}\mu(dx)
 \end{eqnarray}
and
 \begin{eqnarray}\label{3.3}
{\cal B} F(\mu) &=& \frac{1}{2} \int_{\mbb{R}}
\sigma(x)\frac{\delta^2 F(\mu)} {\delta\mu(x)^2} \mu(dx) -
\int_{\mbb{R}} b(x)\frac{\delta F(\mu)} {\delta\mu(x)}
\mu(dx)  \nonumber \\
& &\qquad + \int_{\mbb{R}} \frac{\delta F(\mu)} {\delta\mu(x)} m
(dx)
 \end{eqnarray}
if the right hand sides are meaningful. We shall also prove that
$\{X_t: t\ge0\}$ solves a martingale problem associated with
${\cal L} := {\cal A} + {\cal B}$. It is easily seen that formally
${\cal A} = \lim_{\theta\to0} {\cal A}_\theta$ and ${\cal B} =
\lim_{\theta\to0} {\cal B}_\theta$. Heuristically, $\{X_t:
t\ge0\}$ arises as the high density limit of the immigration
branching particle system discussed in the last section. In
particular, if $F_{f,\{\phi_i\}} (\mu) = f(\<\phi_1,\mu\>, \cdots,
\<\phi_n,\mu\>)$ for $f\in C^2_0(\mbb{R}^n)$ and $\{\phi_i\}
\subset C^2_b(\mbb{R})$, then
 \begin{eqnarray}\label{3.4}
{\cal A} F_{f,\{\phi_i\}}(\mu) &=& \frac{1}{2}\sum_{i,j=1}^n
f_{ij}^{\prime\prime}(\<\phi_1,
\mu\>,\cdots,\<\phi_n,\mu\>)\int_{\mbb{R}^2}\rho(x-y)
\phi_i^\prime(x)\phi_j^\prime(y)\mu(dx)\mu(dy)  \nonumber \\
& &
+\,\frac{1}{2}\sum_{i=1}^n f_i^\prime(\<\phi_1,\mu\>,\cdots,
\<\phi_n,\mu\>)\<a\phi_i^{\prime\prime},\mu\>
 \end{eqnarray}
and
 \begin{eqnarray}\label{3.5}
{\cal B} F_{f,\{\phi_i\}}(\mu) &=& \frac{1}{2}\int_{\mbb{R}}
\sigma(x)\bigg[\sum_{i,j=1}^n
f^{\prime\prime}_{ij}(\<\phi_1,\mu\>,\cdots,
\<\phi_n,\mu\>)\phi_i(x)\phi_j(x)\bigg] \mu(dx)  \nonumber \\
& & - \int_{\mbb{R}} b(x)\bigg[\sum_{i=1}^n f^\prime_i
(\<\phi_1,\mu\>,\cdots,\<\phi_n,\mu\>)\phi_i(x)\bigg]
\mu(dx)  \nonumber \\
& & + \int_{\mbb{R}} \bigg[\sum_{i=1}^n f^\prime_i
(\<\phi_1,\mu\>,\cdots,\<\phi_n,\mu\>)\phi_i(x)\bigg] m(dx).
 \end{eqnarray}
Let ${\cal D}_1({\cal L})$ denote the collection of all functions
$F_{f,\{\phi_i\}}$ with $f\in C^2_0(\mbb{R}^n)$ and $\{\phi_i\}
\subset C^2_b(\mbb{R})$.

We shall obtain (\ref{3.1}) as the limit of a sequence of
equations of immigration branching particle systems. Let $(c, h,
\gamma_k, p^{(k)}, \theta_km, \theta_k^{-1})$ be a sequence of
parameters such that $\theta_k \to \infty$ as $k\to \infty$. Let
$q_k$ and $\sigma_k$ be defined by (\ref{2.5}) and (\ref{2.16}) in
terms of $(\gamma_k, p^{(k)}, \theta_k)$. We assume that
$\{X^{(k)}_t: t\ge0\}$ is a immigration particle system which
satisfies
 \begin{eqnarray}\label{3.6}
\<\phi,X_t^{(k)}\> &=& \<\phi,X_0^{(k)}\> + t\<\phi,m\> +
\int_{(0,t]}\int_{\mbb{R}}\theta_k^{-1}\phi(x)\tilde
N^{(k)}(ds,dx)  \nonumber \\
& & + \int_0^t\<\gamma_k(q_k-1)\phi,X_s^{(k)}\> ds
+ Z_t^{(k)}(\phi) + U_t^{(k)}(\phi^\prime) \nonumber  \\
& & + \int_0^t\int_{\mbb{R}}\<h(y-\cdot)\phi^\prime,X_s^{(k)}\>
W^{(k)}(ds,dy) + \frac{1}{2}\int_0^t \<a\phi^{\prime\prime},
X_s^{(k)}\> ds,
 \end{eqnarray}
where $(N^{(k)},Z^{(k)},M^{(k)},W^{(k)})$ are as in (\ref{2.20})
with parameters $(c, h, \gamma_k, p^{(k)}, \theta_km, \theta_k
^{-1})$. We assume that the $X^{(k)}_0$ are deterministic and
$X^{(k)}_0 \to \mu$ as $k\to \infty$.

\blemma\label{l3.1} Suppose that $B_1 := \sup_{k\ge1}
\|\gamma_k(q_k-1)\| < \infty$ and $B_2 :=\sup_{k\ge1}
\|\theta_k^{-1} \gamma_k \sigma_k\| < \infty$. Then for any
$\phi\in C^2_b(\mbb{R})$, the sequence $\{(\<\phi,X_t^{(k)}\>)
_{t\ge0}, k=1,2,\cdots\}$ is tight in the Skorokhod space
$D([0,\infty), \mbb{R})$. \elemma

\proof Suppose that $\{\tau_k\}$ is a bounded sequence of stopping
times. Let
 \begin{eqnarray*}
V_t^{(k)}(\phi^\prime) = \int_0^t\int_{\mbb{R}}\<h(y-\cdot)
\phi^\prime,X_s^{(k)}\> W^{(k)}(ds,dy)
 \end{eqnarray*}
and
 \begin{eqnarray*}
Y_t^{(k)}(\phi) = \int_0^t\<\gamma_k(q_k-1)\phi,X_s^{(k)}\> ds.
 \end{eqnarray*}
It is easily seen that
 \begin{eqnarray*}
\mbf{E}\{|V_{\tau_k+t}^{(k)}(\phi^\prime) -
V_{\tau_k}^{(k)}(\phi^\prime)|^2\}
 &=&
\mbf{E}\bigg\{\int_0^tds\int_{\mbb{R}}\<h(y-\cdot)
\phi^\prime,X_{\tau_k+s}^{(k)}\>^2 dy\bigg\}  \\
 &=&
\mbf{E}\bigg\{\int_0^tds\int_{\mbb{R}^2}\rho(x-z)\phi^\prime(x)
\phi^\prime(z) X_{\tau_k+s}^{(k)}(dx)X_{\tau_k+s}^{(k)}(dz)\bigg\}  \\
 &\le&
\|\rho\|\int_0^t\mbf{E}\{\<\phi^\prime,X_{\tau_k+s}^{(k)}\>^2\}ds
 \end{eqnarray*}
and
 \begin{eqnarray*}
\mbf{E}\{|Y_{\tau_k+t}^{(k)}(\phi) - Y_{\tau_k}^{(k)}(\phi)|^2\}
\le B_1^2t\int_0^t\mbf{E}\{\<\phi,X_{\tau_k+s}^{(k)}\>^2\} ds.
 \end{eqnarray*}
The remaining terms on the right hand side of (\ref{3.6}) can be
estimated by similar calculations. Combining those estimates and
Proposition~\ref{p2.2} we get
 \begin{eqnarray*}
\sup_{0\le t\le T}\sup_{k\ge1} \mbf{E}\{\<\phi,X_{t}^{(k)}\>^2\} <
\infty
 \end{eqnarray*}
and
 \begin{eqnarray*}
\sup_{k\ge1} \mbf{E}\{|\<\phi,X_{\tau_k+t}^{(k)}\> - \<\phi,
X_{\tau_k}^{(k)}\>|^2\} \to 0
 \end{eqnarray*}
as $t\to 0$. Then the sequence $\{(\<\phi,X_t^{(k)}\>)_{t\ge0},
k=1,2,\cdots\}$ is tight in $D([0,\infty),\mbb{R})$; see
\cite{A78}. \qed

\blemma\label{l3.2} Suppose that $\gamma_k(1-q_k(\cdot)) \to
b(\cdot)$ and $\theta_k^{-1}\gamma_k \sigma_k(\cdot) \to
\sigma(\cdot)$ uniformly for $b\in C_b(\mbb{R})$ and $\sigma\in
C_b(\mbb{R})^+$. Then the sequence $\{X_t^{(k)}: t\ge0,
k=1,2,\cdots\}$ is tight in $D([0,\infty), M(\mbb{R}))$. Moreover,
the limit process $\{X_t: t\ge0\}$ of any subsequence of
$\{X_t^{(k)}: t\ge0, k = 1,2, \cdots\}$ is a.s.\ continuous and
solves the $({\cal L},{\cal D}_1({\cal L}))$-martingale problem,
that is, for each $F \in {\cal D}_1 ({\cal L}_\theta)$,
 \begin{eqnarray}\label{3.7}
F(X_t) - F(X_0) - \int_0^t {\cal L} F(X_s)ds, \qquad t\ge0,
 \end{eqnarray}
is a martingale. \elemma

\proof By Lemma~\ref{l3.1} and a result of \cite{RC86}, the
sequence of processes $\{X_t^{(k)}: t\ge0, k=1,2,\cdots\}$ is
tight in $D([0,\infty), M(\bar\mbb{R}))$. We write $\phi \in C^2_b
(\bar\mbb{R})$ if $\phi \in C^2_b(\mbb{R})$ and its derivatives up
to the second degree can be extended continuously to
$\bar\mbb{R}$. If $\{\phi_i\} \subset C^2 (\bar \mbb{R})$, we can
extend $F_{f,\{\phi_i\}}$, ${\cal A} F_{f,\{\phi_i\}}$ and ${\cal
B} F_{f,\{\phi_i\}}$ continuously to $M(\bar \mbb{R})$. Let $\bar
F_{f,\{\phi_i\}}$, $\bar{\cal A} \bar F_{f,\{\phi_i\}}$ and
$\bar{\cal B} \bar F_{f,\{\phi_i\}}$ denote respectively those
extensions. Let $({\cal A}_k,{\cal B}_k)$ and $(\bar{\cal A}_k,
\bar{\cal B}_k)$ denote the corresponding operators associated
with $\{X_t^{(k)}: t\ge0\}$. Clearly, if $\mu_k\in
M_k(\bar\mbb{R})$ and $\mu_k\to\mu$, then $\bar{\cal A}_k \bar
F_{f,\{\phi_i\}}(\mu_k) \to \bar{\cal A} \bar
F_{f,\{\phi_i\}}(\mu)$. By Taylor's expansion,
 \begin{eqnarray*}
 & & \bar{\cal B}_k \bar F_{f,\{\phi_i\}}(\mu_k)  \\
 &=&
\sum_{j=0}^\infty \int_{\mbb{R}} \theta_k \gamma_k p_j(x)
\big[f(\<\phi_1,\mu_k\>+(j-1)\theta_k^{-1}\phi_1(x),\cdots,
\<\phi_n,\mu_k\>+(j-1)\theta_k^{-1}\phi_n(x))  \nonumber \\
& &\hskip2cm - f(\<\phi_1,\mu_k\>,\cdots,
\<\phi_n,\mu_k\>)\big] \mu_k(dx)  \nonumber \\
& & + \int_{\mbb{R}} \theta_k
\big[f(\<\phi_1,\mu_k\>+\theta_k^{-1}\phi_1(x),\cdots,
\<\phi_n,\mu_k\>+\theta_k^{-1}\phi_n(x))  \nonumber \\
& &\hskip2cm - f(\<\phi_1,\mu_k\>,\cdots, \<\phi_n,\mu_k\>)\big] m
(dx)  \\
 &=&
\int_{\mbb{R}} \gamma_k(q_k(x)-1)\bigg[\sum_{i=1}^n f_i^\prime
(\<\phi_1,\mu_k\>,\cdots,\<\phi_n,\mu_k\>)\phi_i(x)\bigg]
\mu_k(dx)  \nonumber \\
& & + \int_{\mbb{R}} \frac{\gamma_k\sigma_k(x)}{2\theta_k}
\bigg[\sum_{i,j=1}^n f_{ij}^{\prime\prime}(\<\phi_1,\mu_k\>
+\eta_k\phi_1(x),\cdots,\<\phi_n,\mu_k\>+\eta_k\phi_n(x))
\phi_i(x)\phi_j(x)\bigg]\mu_k(dx)  \nonumber \\
& & + \int_{\mbb{R}} \sum_{i=1}^n
\bigg[f^\prime_i(\<\phi_1,\mu_k\>+\zeta_k\phi_1(x),\cdots,
\<\phi_n,\mu_k\>+\zeta_k\phi_n(x))\phi_i(x)\bigg] m (dx),
 \end{eqnarray*}
where $0<\eta_k, \zeta_k<\theta_k^{-1}$. Then $\bar{\cal B}_k \bar
F_{f,\{\phi_i\}}(\mu_k) \to \bar{\cal B} \bar
F_{f,\{\phi_i\}}(\mu)$ under the assumption. Let $\{X_t: t\ge0\}$
be the limit of any subsequence of $\{X_t^{(k)}: t\ge0, k = 1,2,
\cdots\}$. As in the proof of Lemma~4.2 of Dawson \textit{et al}
\cite{DLW01} one can show that
 \begin{eqnarray*}
\bar F_{f,\{\phi_i\}}(X_t) - \bar F_{f,\{\phi_i\}}(X_0) - \int_0^t
\bar{\cal L} \bar F_{f,\{\phi_i\}}(X_s)ds
 \end{eqnarray*}
is a martingale, where $\bar{\cal L} = \bar{\cal A} + \bar{\cal
B}$. As in \cite{W98}, it is not hard to check that the ``gradient
squared'' operator associated with $\bar{\cal L}$ satisfies the
derivation property of \cite{BE85}. Then $\{X_t: t\ge0\}$ is
actually almost surely continuous as an $M(\bar \mbb{R})$-valued
process. By a modification of the proof of Theorem~4.1 of
\cite{DLW01} one can show that $\{X_t: t\ge0\}$ is almost surely
supported by $\mbb{R}$. Thus $\{X_t^{(k)}: t\ge0, k=1,2,\cdots\}$
is tight in $D([0,\infty), M(\mbb{R}))$ and $\{X_t: t\ge0\}$ is
a.s.\ continuous as an $M(\mbb{R})$-valued process. \qed

\blemma\label{l3.3} If $\{X_t: t\ge0\}$ is the continuous solution
of the $({\cal L},{\cal D}_1({\cal L}))$-martingale problem, then
for each integer $n\ge 1$ there is a locally bounded function
$C_n$ on $\mbb{R}_+^3$ such that
 \begin{eqnarray}\label{3.8}
\mbf{E}\{\mbox{$\sup_{0\le s\le t}\<1,X_s\>^n$}\} \le
C_n(\|b\|,\|\sigma\|,t)(1+\<1,\mu\>^n + \<1,m\>^n), \quad t\ge 0.
 \end{eqnarray}
\elemma

\proof If $\{X_t: t\ge0\}$ is the continuous solution of the
$({\cal L},{\cal D}_1({\cal L}))$-martingale problem, then
 \begin{eqnarray}\label{3.9}
Z_t(1) := \<1,X_t\> - \<1,\mu\> - t\<1,m\>
+ \int_0^t\<b,X_s\> ds
 \end{eqnarray}
is a continuous local martingale with quadratic variation process
 \begin{eqnarray}\label{3.10}
\<Z(1)\>_t
=
\int_0^t \<\sigma,X_s\> ds.
 \end{eqnarray}
For $l >0$ let $\tau_l = \inf\{s\ge0: \<1,X_s\> \ge l\}$. The
inequalities for $n=1$ and $n=2$ can be proved as in the proof of
Proposition~\ref{p2.2}. Now the Burkholder-Davis-Gundy inequality
implies that
 \begin{eqnarray*}
\mbf{E}\{\mbox{$\sup_{0\le s\le t}\<1,X_{s\land\tau_l}\>^{2n}$}\}
&\le& C_n\bigg[\<1,\mu\>^{2n} + t^{2n}\<1,m\>^{2n} +
\mbf{E}\bigg\{\bigg(\int_0^{t\land\tau_l}
\<|b|,X_s\>ds\bigg)^{2n}\bigg\}  \\
 & &
+\,\mbf{E}\bigg\{\bigg(\int_0^{t\land\tau_l} \<\sigma,X_s\>
ds\bigg)^n\bigg\}\bigg] \\
 &\le&
C_n\bigg[\<1,\mu\>^{2n} + t^{2n}\<1,m\>^{2n}
+ \theta^{-n} t^n\<1,m\>^n \nonumber \\
 & &
+\,\|b\|^{2n}t^{2n-1}\int_0^t \mbf{E}\{\mbox{$\sup_{0\le r\le
s}\<1,X_{r\land\tau_l}\>^{2n}$}\}
ds\bigg] \nonumber \\
 & &
+\,\|\sigma\|^n t^{n-1}\int_0^t\mbf{E}\{\<1,X_s\>^n\} ds,
 \end{eqnarray*}
where $C_n\ge0$ is a universal constant. By using the above
estimate and Gronwall's inequality inductively, we get some
estimates for $\mbf{E}\{\sup_{0\le s\le t} \<1,$ $X_{t \land
\tau_l}\>^n\}$. Then we obtain the inequalities for
$\mbf{E}\{\sup_{0\le s\le t} \<1,X_t\>^n\}$ by Fatou's lemma. \qed

\blemma\label{l3.4} Suppose there are constants $d_0>0$ and
$\delta>1/2$ such that $h(x)\le d_0(1+|x|)^{-\delta}$ for all
$x\in \mbb{R}$. If $\gamma_k (1-q_k(\cdot)) \to b(\cdot)$ and
$\theta_k^{-1} \gamma_k \sigma_k (\cdot) \to \sigma(\cdot)$
uniformly for $b\in C_b(\mbb{R})$ and $\sigma \in C_b(\mbb{R})^+$,
then the limit process $\{X_t: t\ge0\}$ of any subsequence of
$\{X_t^{(k)}: t\ge0, k= 1,2, \cdots\}$ is a weak solution of
(\ref{3.1}). \elemma

\proof By the proof of Lemma~\ref{l3.1} and the results of
\cite{M83, RC86}, $\{(X_t^{(k)}, U_t^{(k)}, W_t^{(k)}, Z_t^{(k)}):
t\ge0, k=1,2,\cdots\}$ is a tight  sequence in $D([0,\infty)$,
$M(\bar\mbb{R})\times {\cal S}^\prime (\mbb{R})^3)$. By passing to
a subsequence, we simply assume that $\{(X_t^{(k)}, U_t^{(k)},
W_t^{(k)}, Z_t^{(k)}): t\ge0\}$ converges in distribution to some
process $\{(X_t, U_t, W_t, Z_t): t\ge0\}$. By Lemma~\ref{l3.2},
$\{X_t: t\ge0\}$ is a.s.\ continuous and solves the $({\cal
L},{\cal D}_1({\cal L}))$-martingale problem. Considering the
Skorokhod representation, we assume $\{(X_t^{(k)}, U_t^{(k)},
W_t^{(k)}, Z_t^{(k)}): t\ge0\}$ converges almost surely to the
process $\{(X_t, U_t, W_t, Z_t): t\ge0\}$ in the topology of
$D([0,\infty), M(\bar\mbb{R})\times {\cal S}^\prime(\mbb{R})^3)$.
Since each $\{W_t^{(k)}: t\ge0\}$ is a time-space white noise, so
is $\{W_t: t\ge0\}$. In view of (\ref{2.15}), we have a.s.\
$U_t(\phi)=0$ for all $t\ge0$ and $\phi \in {\cal S}(\mbb{R})$.
Then the theorem follows once it is proved that $\{(X_t, W_t,
Z_t): t\ge0\}$ satisfies (\ref{3.1}). Clearly, it is sufficient to
prove this for $\phi \in {\cal S}(\mbb{R})$ with compact support
$\mbox{supp} (\phi)$. Let $Y_t(y) = \<h(y-\cdot) \phi^\prime,
X_t\>$ and $Y_t^{(k)}(y) = \<h(y-\cdot) \phi^\prime, X_t^{(k)}\>$.
For $l >0$ let $\tau_l = \inf\{s\ge0: \<1,X_s^{(k)}\> \ge l$ for
some $k\ge1\}$. Since the weak convergence of measures can be
induced by the (Vasershtein) metric defined in \cite[p.150]{EK86},
it is easy to show that $\{Y_t^{(k)}1_{\{t<\tau_l\}}: t\ge0\}$
converges to $\{Y_t1_{\{t<\tau_l\}}: t\ge0\}$ in $D([0,\infty),
C_0(\mbb{R}))$, where $C_0(\mbb{R})$ is furnished with the uniform
norm. By \cite[Theorem~2.1]{C95}, for $\psi\in{\cal S}(\mbb{R})$
we have almost surely
 \begin{eqnarray}\label{3.11}
\lim_{k\to\infty}\int_0^t\int_{\mbb{R}}\psi(y)Y_s^{(k)}(y)
1_{\{s<\tau_l\}} W^{(k)}(ds,dy) = \int_0^t\int_{\mbb{R}}
\psi(y)Y_s(y)1_{\{s<\tau_l\}} W(ds,dy).
 \end{eqnarray}
Let $\alpha = \sup\{|x|, x\in \mbox{supp} (\phi)\}$. We have
 \begin{eqnarray*}
\sup_{|z|\le \alpha}|h(y-z)| \le d(y) := d_0[1_{\{|y|\le \alpha\}}
+ 1_{\{|y|>\alpha\}}(1+|y|-\alpha)^{-\delta}],
 \end{eqnarray*}
and hence
 \begin{eqnarray*}
|Y_t(y)| \le \<|\phi^\prime|, X_t\>d(y) \quad\mbox{and}\quad
|Y_t^{(k)}(y)| \le \<|\phi^\prime|, X_t^{(k)}\> d(y).
 \end{eqnarray*}
By the Burkholder-Davis-Gundy inequality,
 \begin{eqnarray}\label{3.12}
&&\mbf{E}\bigg\{\bigg(\int_0^t\int_{\mbb{R}}\psi(y) Y_s^{(k)}(y)
1_{\{s<\tau_l\}}W^{(k)}(ds,dy)\bigg)^4\bigg\} \nonumber \\
&&\hskip2cm\le
\mbox{const}\cdot\mbf{E}\bigg\{\bigg(\int_0^t\int_{\mbb{R}}
\psi(y)^2Y_s^{(k)}(y)^2
1_{\{s<\tau_l\}}ds dy\bigg)^2\bigg\} \nonumber \\
&&\hskip2cm\le
\mbox{const}\cdot l^4\|\phi^\prime\|^4\<\psi^2d^2,\lambda\>^2t^2,
 \end{eqnarray}
where $\lambda$ denotes the Lebesgue measure on $\mbb{R}$. Since
the right hand side of (\ref{3.12}) is independent of $k\ge1$, the
convergence of (\ref{3.11}) also holds in the $L^2$-sense. For
each $\epsilon>0$, it is not hard to choose $\psi\in {\cal
S}(\mbb{R})$ so that
 \begin{eqnarray}\label{3.13}
&&\mbf{E}\bigg\{\bigg(\int_0^t\int_{\mbb{R}}(1-\psi(y))
Y_s^{(k)}(y)
1_{\{s<\tau_l\}}W^{(k)}(ds,dy)\bigg)^2\bigg\} \nonumber \\
&&\hskip2cm \le \mbox{const}\cdot
l^2\|\phi^\prime\|^2\<|1-\psi|^2d^2,\lambda\>t \le \epsilon.
 \end{eqnarray}
The same estimate is available with $Y^{(k)}$ and $W^{(k)}$
replaced respectively by $Y$ and $W$. Clearly, (\ref{3.11}) and
(\ref{3.13}) imply that
 \begin{eqnarray}\label{3.14}
\lim_{k\to\infty}\int_0^t\int_{\mbb{R}}Y_s^{(k)}(y)
1_{\{s<\tau_l\}}W^{(k)}(ds,dy) = \int_0^t\int_{\mbb{R}}Y_s(y)
1_{\{s<\tau_l\}}W(ds,dy)
 \end{eqnarray}
in the $L^2$-sense. Passing to a suitable subsequence we get the
almost sure convergence for (\ref{3.14}). Now letting $k\to
\infty$ in (\ref{3.6}) we get
 \begin{eqnarray*}
\<\phi,X_{t\land\tau_l}\> &=& \<\phi,\mu\> + (t\land\tau_l)
\<\phi,m\> + \frac{1}{2}\int_0^{t\land\tau_l}
\<a\phi^{\prime\prime}, X_s\> ds
- \int_0^{t\land\tau_l}\<b\phi,X_s\> ds \nonumber \\
& & + \int_0^{t\land\tau_l}\int_{\mbb{R}}\phi(y) Z(ds,dy) +
\int_0^{t\land\tau_l}\int_{\mbb{R}}\<h(y-\cdot) \phi^\prime,X_s\>
W(ds,dy),
 \end{eqnarray*}
from which (\ref{3.1}) follows. The extensions from $\phi\in {\cal
S}(\mbb{R})$ to $\phi\in C^2_b(\mbb{R})$ is immediate. \qed

\btheorem\label{t3.1} Suppose there are constants $d_0>0$ and
$\delta>1/2$ such that $h(x)\le d_0(1+|x|)^{-\delta}$ for all
$x\in \mbb{R}$. Then the stochastic equation (\ref{3.1}) has a
continuous weak solution $\{X_t: t\ge0\}$. Moreover, $\{X_t:
t\ge0\}$ also solves the $({\cal L}, {\cal D}_1({\cal
L}))$-martingale problem. \etheorem

\proof Given $b\in C_b(\mbb{R})$ and $\sigma\in C_b(\mbb{R})^+$,
we set $\theta_k=k$, $\gamma_k = \sqrt{k}$ and
 \begin{eqnarray*}
p_0^{(k)} = 1 - p_2^{(k)} - p_k^{(k)},
\quad
p_2^{(k)} = \frac{(k-1)^2(1-b/\sqrt{k})-k\sigma_k}{2(k-1)^2-k},
\quad
p_k^{(k)} = \frac{2\sigma_k-1+b/\sqrt{k}}{2(k-1)^2-k},
 \end{eqnarray*}
where $\sigma_k(\cdot) = \sqrt{k}\sigma(\cdot) + 1$. Then the
sequence $(\gamma_k, p^{(k)}, \theta_k)$ satisfies the conditions
of Lemma~\ref{l3.4}. By Lemmas~\ref{l3.2} and \ref{l3.4}, equation
(\ref{3.1}) has a continuous weak solution $\{X_t: t\ge0\}$ which
solves the $({\cal L},{\cal D}_1({\cal L}))$-martingale problem.
\qed


\section{Stochastic log-Laplace equations}

\setcounter{equation}{0}

In this section, we establish the existence and uniqueness of
solution of the stochastic log-Laplace equation (\ref{1.8}). The
techniques here are based on the results of Kurtz and Xiong
\cite{KX99} and have been stimulated by \cite{C04, X04}. Let $(c,
h, \sigma, b, m)$ be given as in the introduction. Suppose that
$W(ds,dx)$ is a time-space white noise. The main objective is to
discuss the non-linear SPDE:
 \begin{eqnarray}\label{4.1}
 \psi_t(x) &=& \phi(x)
+ \int_0^t \bigg[\frac{1}{2}a(x)\partial_x^2\psi_s(x) -
b(x)\psi_s(x) - \frac{1}{2}\sigma(x)\psi_s(x)^2\bigg] ds
\nonumber \\
& & + \int_0^t\int_{\mbb{R}} h(y-x)
\partial_x\psi_s(x) W(ds,dy),
\qquad t\ge 0.
 \end{eqnarray}

Let $\{H_k(\mbb{R}): k=0,\pm 1, \pm 2, \cdots\}$ denote the
Sobolev spaces on $\mbb{R}$. Let ``$\|\cdot\|_0$'' and ``$\<\cdot,
\cdot\>_0$'' denote respectively the norm and the inner product in
$H_0(\mbb{R}) = L^2(\mbb{R})$. For $\phi\in H_k(\mbb{R})$ let
 \begin{eqnarray}\label{4.2}
\|\phi\|_k^2 = \sum_{i=0}^k \|\partial_x^{i}\phi\|_0^2.
 \end{eqnarray}
Following Xiong \cite{X04}, we first consider a smoothed version
of equation (\ref{4.1}). Let $(T_t)_{t\ge0}$ denote the transition
semigroup of a standard Brownian motion. Let $\{h_j:
j=1,2,\cdots\}$ be a complete orthonormal system of
$H_0(\mbb{R})$. Then
 \begin{eqnarray}\label{4.3}
W_j(t) =\int_0^t\int_{\mbb{R}}h_j(y)W(ds,dy), \qquad t\ge 0
 \end{eqnarray}
defines a sequence of independent standard Brownian motions
$\{W_j: j=1,2,\cdots\}$. For $\epsilon>0$ let
 \begin{eqnarray}\label{4.4}
W^\epsilon(dt,dx) = \sum_{j=1}^{[1/\epsilon]} h_j(x)W_j(dt)dx,
\qquad s\ge0, y\in\mbb{R}.
 \end{eqnarray}
For $\phi\in H_0(\mbb{R})$ we set $d_\epsilon(\phi) =
(\|T_\epsilon \phi\| \land \epsilon^{-1}) \|T_\epsilon
\phi\|^{-1}$. By the general results of \cite[Theorem~3.5]{KX99}
and \cite[p.133]{R90}, for any $\phi\in H_1(\mbb{R}) \cap
C_b(\mbb{R})^+$ there is a pathwise unique $H_2(\mbb{R})$-valued
solution $\{\psi_t^\epsilon: t\ge0\}$ of the equation
 \begin{eqnarray}\label{4.5}
\psi_t^\epsilon(x) &=& T_{\epsilon}\phi(x) + \int_0^t
\bigg[\frac{1}{2}a(x)\partial_x^2\psi_s^\epsilon(x) -
b(x)\psi_s^\epsilon(x) - \frac{1}{2}\sigma(x) \psi_s^\epsilon(x)
d_\epsilon(\psi^\epsilon_s) T_\epsilon \psi_s^\epsilon(x)\bigg] ds
\nonumber \\
& & + \int_0^t\int_{\mbb{R}} h(y-x)\partial_x\psi_s^\epsilon(x)
W^\epsilon(ds,dy), \qquad t\ge 0.
 \end{eqnarray}

\blemma\label{l4.1} The solution $\{\psi_t ^\epsilon: t\ge0\}$ of
(\ref{4.5}) is non-negative and satisfies a.s.\
$\|\psi_t^\epsilon\|_{\mbox{\rm ess}} \le e^{-b_0t}
\|\phi\|_{\mbox{\rm ess}}$ for all $t\ge0$, where $b_0 = \inf_x
b(x)$ and $\|\cdot\|_{\mbox{\rm ess}}$ denote the essential
supremum norm. \elemma

\proof Indeed, for any non-negative and non-trivial function
$\phi\in H_0(\mbb{R})$, the solution of (\ref{4.5}) can be
obtained in the following way. Let $\{B_i(t)\}$ be a sequence of
independent Brownian motions which are also independent of the
white noise $\{W(ds,dy)\}$. As in \cite[Theorems~2.1 and
2.2]{KX99}, one can show that there is a pathwise unique solution
$\psi_t^\epsilon(x)$ of the stochastic system
 \begin{eqnarray}\label{4.6}
\xi_i(t) - \xi_i(0)
&=&
\int_0^t c(\xi_i(s))dB_i(s)
+ 2\int_0^t c(\xi_i(s))c^\prime(\xi_i(s))ds \hskip1cm \nonumber \\
& & -\int_0^t\int_{\mbb{R}} h(y-\xi_i(s)) W^\epsilon(ds,dy),
 \end{eqnarray}
 \begin{eqnarray}\label{4.7}
m_i(t) - m_i(0) &=& \int_0^t \bigg[\frac{1}{2}a^{\prime\prime}
(\xi_i(s)) - b(\xi_i(s))\bigg]m_i(s)ds  \nonumber \\
& & -
\frac{1}{2}\int_0^t\sigma(\xi_i(s))d_\epsilon(\psi_s^\epsilon)
T_\epsilon\psi_s^\epsilon(\xi_i(s)) m_i(s)ds  \nonumber \\
& & - \int_0^t\int_{\mbb{R}} h^\prime(y-\xi_i(s)) m_i(s)
W^\epsilon(ds,dy),
 \end{eqnarray}
and
 \begin{eqnarray}\label{4.8}
\psi_t^\epsilon(x)dx = \lim_{n\to\infty}\frac{1}{n}\sum_{i=1}^n
m_i(t) \delta_{\xi_i(t)}(dx), \qquad t\ge0, x\in \mbb{R},
 \end{eqnarray}
where $\{(m_i(0),\xi_i(0)): i=1,2,\cdots\}$ is a sequence of
exchangeable random variables on $[0,\infty)\times \mbb{R}$ which
are independent of $\{B_i(t)\}$ and $\{W(ds,dy)\}$ and satisfy
 $$
\lim_{n\to\infty} n^{-1}\sum_{i=1}^n m_i(0) \delta_{\xi_i(0)}(dx)
= T_{\epsilon}\phi(x)dx.
 $$
By the arguments of \cite[Theorems~3.1-3.5]{KX99}, it can be
proved that $\psi_t^\epsilon (x)$ is also the pathwise unique
solution of (\ref{4.5}). By a duality argument similar to the
proof of \cite[Lemma~2.2]{X04} we get $\|\psi_t^\epsilon\|
_{\mbox{ess}} \le e^{-b_0t} \|\phi\|_{\mbox{ess}}$. \qed

\blemma\label{l4.2} There is a locally bounded function $K(\cdot)$
on $[0,\infty)$ such that
 \begin{eqnarray}\label{4.9}
\mbf{E}\bigg\{\sup_{0\le r\le t} \|\psi_r^\epsilon\|_0^4\bigg\}
\le K(t), \qquad t\ge0.
 \end{eqnarray}
\elemma

\proof Although the arguments are similar to those of \cite{X04},
we shall give the detailed proof for the convenience of the
reader. For any $f\in C^\infty(\mbb{R})$ with compact support,
 \begin{eqnarray*}
\<\psi_t^\epsilon,f\>_0
 &=&
\<T_{\epsilon}\phi,f\>_0 + \int_0^t
\bigg[\frac{1}{2}\<a\partial_x^2 \psi_s^\epsilon,f\>_0 -
\<b\psi_s^\epsilon,f\>_0 - \frac{1}{2}\<\sigma\psi_s^\epsilon
d_\epsilon(\psi_s^\epsilon)
T_\epsilon\psi_s^\epsilon,f\>_0\bigg]ds \nonumber \\
 & &
+ \int_0^t\int_{\mbb{R}} \<h(y-\cdot)\partial_x \psi_s^\epsilon,
f\>_0 W^\epsilon (ds,dy).
 \end{eqnarray*}
By It\^o's formula,
 \begin{eqnarray*}
\<\psi_t^\epsilon,f\>_0^2 &=& \<T_{\epsilon}\phi,f\>_0^2 +
\int_0^t \<\psi_s^\epsilon,f\>_0 \<a\partial_x^2 \psi_s^\epsilon -
2b\psi_s^\epsilon - \sigma\psi_s^\epsilon
d_\epsilon(\psi_s^\epsilon)T_\epsilon\psi_s^\epsilon,
f\>_0 ds \nonumber \\
& & +\, 2\int_0^t\int_{\mbb{R}}\<\psi_s^\epsilon,f\>_0
\<h(y-\cdot)\partial_x\psi_s^\epsilon,f\>_0 W^\epsilon(ds,dy)
\nonumber \\
& & + \sum_{j=1}^{[1/\epsilon]}\int_0^t\bigg[\int_{\mbb{R}}
h_j(y)\<h(y-\cdot)\partial_x\psi_s^\epsilon,f\>_0 dy\bigg]^2 ds.
 \end{eqnarray*}
Then we may add $f$ over in a complete orthonormal system of
$H_0(\mbb{R})$ to get
 \begin{eqnarray}\label{4.10}
\|\psi_t^\epsilon\|_0^2 &=& \|T_{\epsilon}\phi\|_0^2 + \int_0^t
\<a\partial_x^2\psi_s^\epsilon - 2b\psi_s^\epsilon -
\sigma\psi_s^\epsilon d_\epsilon(\psi_s^\epsilon) T_\epsilon
\psi_s^\epsilon,\psi_s^\epsilon\>_0 ds \nonumber \\
 & &
+\, 2\int_0^t\int_{\mbb{R}}\<h(y-\cdot)\partial_x \psi_s^\epsilon,
\psi_s^\epsilon\>_0 W^\epsilon(ds,dy) \nonumber \\
 & &
+ \sum_{j=1}^{[1/\epsilon]}\int_0^t ds \int_{\mbb{R}} \bigg[
\int_{\mbb{R}} h_j(y) h(y-z)\partial_x \psi_s^\epsilon(z)
dy\bigg]^2dz \nonumber\\
 &\le&
\|T_{\epsilon}\phi\|_0^2 + \int_0^t\<c^2\partial_x^2
\psi_s^\epsilon, \psi_s^\epsilon\>_0 ds + \int_0^t
\<\rho(0)\partial_x^2\psi_s^\epsilon,
\psi_s^\epsilon\>_0 ds  \nonumber \\
 & &
+ \int_0^t \< - 2b\psi_s^\epsilon - \sigma\psi_s^\epsilon
d_\epsilon(\psi_s^\epsilon) T_\epsilon
\psi_s^\epsilon,\psi_s^\epsilon\>_0 ds \nonumber \\
 & &
+\, 2\int_0^t\int_{\mbb{R}} \<h(y-\cdot)\partial_x
\psi_s^\epsilon,\psi_s^\epsilon\>_0 W^\epsilon(ds,dy) \nonumber \\
 & &
+ \int_0^t ds \int_{\mbb{R}} \bigg[ \int_{\mbb{R}}
h(y-z)^2(\partial_x \psi_s^\epsilon(x))^2 dy\bigg] dx.
 \end{eqnarray}
Note that the third and the last terms on the right hand side
cancel out. Since $\psi^{\epsilon}_s\in H_2(\mbb{R})$, there
exists a sequence $f_n\in C^{\infty}_0(\mbb{R})$ such that
$f_n\to\psi^{\epsilon}_s$ in $H_2(\mbb{R})$. By the assumption,
both $c^2$ and $(c^2)''$ are bounded. Then there is a constant
$K\ge0$ such that
 \begin{eqnarray*}
\<c^2f''_n,f_n\>  = \<(c^2)'',f_n^2\>/2  - \<c^2,(f'_n)^2\> \le
K\|f_n\|^2_0.
 \end{eqnarray*}
Taking $n\to\infty$ we have
 \begin{eqnarray}\label{4.11}
\<c^2\partial_x^2\psi_s^\epsilon, \psi_s^\epsilon\>_0\le
K\|\psi_s^\epsilon\|^2_0.
 \end{eqnarray}
By Lemma~\ref{l4.1}, it is easy to find a locally bounded
non-negative function $K(\cdot)$ such that
 \begin{eqnarray*}
\< - 2b\psi_s^\epsilon - \sigma\psi_s^\epsilon
d_\epsilon(\psi_s^\epsilon) T_\epsilon \psi_s^\epsilon,
\psi_s^\epsilon\>_0 \le K(s)\|\psi_s^\epsilon\|^2_0.
 \end{eqnarray*}
Therefore, we can redesign $K(\cdot)$ suitably and get from
(\ref{4.10}) that
 \begin{eqnarray*}
\|\psi_t^\epsilon\|_0^2
 &\le&
\|\phi\|_0^2 + K(t)\int_0^t \|\psi_s^\epsilon\|^2_0ds +
2\int_0^t\int_{\mbb{R}} \<h(y-\cdot)\partial_x\psi_s^\epsilon,
\psi_s^\epsilon\>_0 W^\epsilon(ds,dy).
 \end{eqnarray*}
By Schwarz' and Burkholder's inequalities we can redesign
$K(\cdot)$ again to get
 \begin{eqnarray}\label{4.12}
\mbf{E}\bigg\{\sup_{0\le r\le t}\|\psi_r^\epsilon\|_0^4\bigg\}
 &\le&
3\|\phi\|_0^4 + K(t)\mbf{E}\bigg\{\int^t_0
\|\psi_s^\epsilon\|^4_0 ds\bigg\} \nonumber\\
 & &
+\, 24\mbf{E}\bigg\{\int_0^t\int_{\mbb{R}} \<h(y-\cdot)
\partial_x\psi_s^\epsilon,\psi_s^\epsilon\>_0^2
dy ds\bigg\} \nonumber\\
 &\le&
3\|\phi\|_0^4 + K(t)\mbf{E}\bigg\{\int^t_0
\|\psi_s^\epsilon\|^4_0ds\bigg\},
 \end{eqnarray}
where the last inequality follows from the same arguments as those
leading to (\ref{4.11}). Using stopping times if necessary, we may assume
that $\mbf{E}\{\|\psi_t^\epsilon\|_0^4\} < \infty$ for each $t\ge 0$. Then we
obtain (\ref{4.9}) by Gronwall's inequality. \qed

\blemma\label{l4.3} There is a locally bounded function $K(\cdot)$
on $[0,\infty)$ such that
 \begin{eqnarray}\label{4.13}
\mbf{E}\bigg\{\sup_{0\le r\le t} \|\psi_r^\epsilon\|_1^4\bigg\}
\le K(t), \qquad t\ge0.
 \end{eqnarray}
\elemma

\proof We shall omit some details since they are similar to those
in the proof of Lemma~\ref{l4.2}. {From} (\ref{4.5}) it follows
that
 \begin{eqnarray*}
\partial_x\psi_t^\epsilon(x)
 &=&
\partial_xT_{\epsilon}\phi(x) + \int_0^t\bigg[\frac{1}{2}
a^\prime(x)\partial_x^2\psi_s^\epsilon(x) + \frac{1}{2}a(x)
\partial_x^3\psi_s^\epsilon(x) - b^\prime(x)\psi_s^\epsilon(x) -
b(x)\partial_x\psi_s^\epsilon(x)  \\
 & &
- \frac{1}{2}\sigma^\prime(x)\psi_s^\epsilon(x)d_\epsilon
(\psi^\epsilon_s) T_\epsilon \psi_s^\epsilon(x) -
\frac{1}{2}\sigma(x)\partial_x \psi_s^\epsilon(x)
d_\epsilon(\psi^\epsilon_s) T_\epsilon\psi_s^\epsilon(x)  \\
 & &
- \frac{1}{2}\sigma(x) \psi_s^\epsilon(x)d_\epsilon
(\psi^\epsilon_s)T_\epsilon\partial_\cdot
\psi_s^\epsilon(x)\bigg] ds   \\
 & &
+ \int_0^t\int_{\mbb{R}} [h(y-x)\partial_x^2\psi_s^\epsilon(x) -
h^\prime(y-x)\partial_x\psi_s^\epsilon(x)] W^\epsilon(ds,dy).
 \end{eqnarray*}
Then we have
 \begin{eqnarray*}
\|\partial_x\psi_t^\epsilon\|_0^2 &=&
\|T_{\epsilon}\partial_x\phi\|_0^2 + \int_0^t
\bigg[\<\partial_x\psi_s^\epsilon,a^\prime\partial_x^2
\psi_s^\epsilon + a\partial_x^3\psi_s^\epsilon\>_0 - 2
\<\partial_x\psi_s^\epsilon, b^\prime\psi_s^\epsilon
+ b \partial_x \psi_s^\epsilon\>_0  \\
 & & - d_\epsilon(\psi^\epsilon_s) \<\partial_x\psi_s^\epsilon,
\sigma^\prime \psi_s T_\epsilon \psi_s^\epsilon + \sigma
\partial_x \psi_s^\epsilon T_\epsilon \psi_s^\epsilon
+ \sigma \psi_s^\epsilon T_\epsilon\partial_x
\psi_s^\epsilon\>_0\bigg] ds \nonumber \\
 & & + 2\int_0^t\int_{\mbb{R}} \<\partial_x\psi_s^\epsilon,
h(y-\cdot)\partial_x^2\psi_s^\epsilon -
h^\prime(y-\cdot)\partial_x
\psi_s^\epsilon\>_0 W^{\epsilon}(ds,dy) \nonumber \\
 & & + \int_0^tds\int_{\mbb{R}} \|h(y-\cdot)\partial_x^2\psi_s^\epsilon
- h^\prime(y-\cdot)\partial_x\psi_s^\epsilon\|_0^2 dy.
 \end{eqnarray*}
As in the proof of the previous lemma, we have that
 \begin{eqnarray}\label{4.14}
 \mbf{E}\bigg\{\sup_{0\le r\le t}\|\partial_x \psi_t^\epsilon
\|_0^4\bigg\} \le4\|\partial_x\phi\|_0^4 + K(t)
\mbf{E}\int^t_0\left(\|\psi_s^\epsilon\|^4_0
+\|\partial_x\psi_s^\epsilon\|^4_0\right)ds.
 \end{eqnarray}
Again, we may assume $\mbf{E}\bigg\{\sup_{0\le r\le t} \|\partial_x
\psi_r^\epsilon \|_0^4\bigg\}<\infty$ for all $t\ge 0$. Then we
obtain (\ref{4.13}) by Gronwall's inequality. \qed

\btheorem\label{t4.1} For any $\phi \in H_1(\mbb{R}) \cap
C_b(\mbb{R})^+$, equation (\ref{4.1}) has a pathwise unique
$H_1(\mbb{R})^+$-valued solution $\{\psi_t: t\ge0\}$. We have
a.s.\ $\|\psi_t\|_{\mbox{\rm ess}} \le e^{-b_0t} \|\phi\|_{\mbox
{\rm ess}}$ for all $t\ge0$. Moreover, there is a locally bounded
function $K(\cdot)$ on $[0,\infty)$ such that
 \begin{eqnarray}\label{4.15}
\mbf{E}\Big\{\sup_{0\le r\le t} \|\psi_r\|_1^4\Big\} \le K(t),
 \end{eqnarray}
and so $\{\psi_t(\cdot): t\ge0\}$ has an $H_1(\mbb{R}) \cap
C_b(\mbb{R})^+$-valued version. \etheorem

\proof Let $z_t(x) = \psi^{\epsilon}_t(x) - \psi^{\eta}_t(x)$. For
any $t\ge0$, by the same arguments leading to (2.12) of \cite{X04}
we have
 \begin{eqnarray}\label{4.16}
\mbf{E}\Big\{\sup_{0\le s\le t}\|z_s\|^4_0\Big\}
 &\le&
K\int^t_0\mbf{E}\{\|z_{r}\|^4_0\}dr + K\mbf{E}
\bigg\{\int^t_0|d_{\epsilon}(\psi^{\epsilon}_r) -
d_{\eta}(\psi^{\eta}_r)|^4dr\bigg\} \nonumber \\
 & &
+\, 3\|\phi\|^4 \mbf{E} \bigg\{\int^t_0\left(\int|T_{\epsilon}
\psi_r^{\epsilon}(x) - T_{\eta}\psi^{\eta}_{r}(x)|^2
dx\right)^2dr\bigg\}  \nonumber \\
 & &
+\, K\mbf{E}\bigg\{\sum^{[1/\epsilon]}_{j=[1/\eta]+1}\int^t_0
\left(\int_{\mbb{R}} \left<h(y-\cdot)\partial_x\psi^{\eta}_s,
z_s\right> h_j(y)dy\right)^2ds\bigg\}.
 \end{eqnarray}
As in Section 2.4 of \cite{X04}, the second and third terms on the
right hand side of (\ref{4.16}) converge to zero as $\epsilon$ and
$\eta\to 0$. On the other hand, the last term is bounded by
 \begin{eqnarray*}
\int^t_0\int_{\mbb{R}}\sum^{[1/\epsilon]}_{j=[1/\eta]+1}
\left(\int_{\mbb{R}} h_j(y)h(y-x)dy\right)^2\mbf{E} \{z_s(x)^2\}
dx \int_{\mbb{R}}\mbf{E}\{(\partial_x\psi^{\eta}_s(x))^2\}dx ds,
 \end{eqnarray*}
which tends to zero $\epsilon$ and $\eta\to 0$. As in Section 2.4
of \cite{X04} we can show that $\psi^{\epsilon}$ is a Cauchy
sequence in $H_0(\mbb{R})$ and its limit $\psi$ is the pathwise
unique solution of (\ref{4.1}). The second assertion follows from
Lemma~\ref{l4.1} and Fatou's lemma. Finally, we obtain
(\ref{4.15}) by Lemma~\ref{l4.3} and Sobolev's result. \qed

Based on Theorem~\ref{t4.1}, let us consider the following more
useful backward SPDE:
 \begin{eqnarray}\label{4.17}
\psi_{r,t}(x) &=& \phi(x) + \int_r^t
\bigg[\frac{1}{2}a(x)\partial_x^2 \psi_{s,t}(x) -
b(x)\psi_{s,t}(x) - \frac{1}{2}\sigma(x) \psi_{s,t}(x)^2\bigg] ds
\nonumber \\
& & + \int_r^t\int_{\mbb{R}} h(y-x)\partial_x\psi_{s,t}(x) \cdot
W(ds,dy), \qquad t\ge r\ge 0,
 \end{eqnarray}
where ``$\cdot$'' denotes the backward stochastic integral.

\btheorem\label{t4.2} For any $\phi \in H_1(\mbb{R}) \cap
C_b(\mbb{R})^+$, the backward equation (\ref{4.17}) has a pathwise
unique $H_1(\mbb{R}) \cap C_b(\mbb{R})^+$-valued solution
$\{\psi_{r,t}: t\ge r\ge0\}$. Further, we have a.s.\ $\|\psi_{r,t}
\| \le e^{-b_0(t-r)} \|\phi\|$ for all $t\ge r\ge0$. \etheorem

\proof For fixed $t>0$, define the white noise
 \begin{eqnarray}\label{4.18}
W_t([0,s]\times B) = - W([t-s,t]\times B), \qquad 0\le s\le t,
B\in {\cal B}(\mbb{R}).
 \end{eqnarray}
By Theorem~\ref{t4.1}, there is a pathwise unique solution
$\{\phi_{r,t}: 0\le r\le t\}$ of the equation
 \begin{eqnarray}\label{4.19}
\phi_{r,t}(x) &=& \phi(x) + \int_0^r
\bigg[\frac{1}{2}a(x)\partial_x^2\phi_{s,t}(x) - b(x)\phi_{s,t}(x)
- \frac{1}{2}\sigma(x)\phi_{s,t}(x)^2\bigg] ds
\nonumber \\
& & + \int_0^r\int_{\mbb{R}} h(y-x)
\partial_x\phi_{s,t}(x) W_t(ds,dy).
 \end{eqnarray}
Setting $\psi_{r,t}(x) := \phi_{t-r,t}(x)$, we have
 \begin{eqnarray*}
\psi_{r,t}(x) &=& \phi(x) + \int_0^{t-r}
\bigg[\frac{1}{2}a(x)\partial_x^2 \psi_{t-s,t}(x) -
b(x)\psi_{t-s,t}(x) - \frac{1}{2}
\sigma(x)\psi_{t-s,t}(x)^2\bigg] ds \nonumber \\
& & + \int_0^{t-r}\int_{\mbb{R}} h(y-x)
\partial_x\psi_{t-s,t}(x) W_t(ds,dy) \nonumber \\
&=&  \phi(x) + \int_r^t \bigg[\frac{1}{2}a(x)\partial_x^2
\psi_{s,t}(x) - b(x)\psi_{s,t}(x) - \frac{1}{2}\sigma(x)
\psi_{s,t}(x)^2\bigg] ds \nonumber \\
& & + \int_r^t\int_{\mbb{R}} h(y-x)\partial_x\psi_{s,t}(x) \cdot
W(ds,dy).
 \end{eqnarray*}
That is, $\{\psi_{r,t}: t\ge r\ge0\}$ solves (\ref{4.17}). The
remaining assertions are immediate by Theorem~\ref{t4.1}. \qed

We may regard the white noise $\{W(ds,dy)\}$ as a random variable
taking values in the Schwartz apace ${\cal S}^\prime([0,\infty)
\times \mbb{R})$. As in the classical situation of
\cite[p.163]{IW89}, the result of Theorem~\ref{t4.2} implies the
existence of a measurable mapping $F: (\phi,w) \mapsto
\psi_{r,t}^w(\phi,\cdot)$ from $(H_1(\mbb{R}) \cap C_b(\mbb{R})^+)
\times {\cal S}^\prime([0,\infty) \times \mbb{R})$ to
$H_1(\mbb{R}) \cap C_b(\mbb{R})^+$ such that
$\psi_{r,t}^W(\phi,\cdot)$ is the pathwise unique solution of
(\ref{4.17}).


\section{Conditional log-Laplace functionals}

\setcounter{equation}{0}

Let $(c, h, \sigma, b, m)$ be given as in the introduction. Let
$\{X_t: t\ge0\}$ be a continuous solution of the SPDE:
 \begin{eqnarray}\label{5.1}
\<\phi,X_t\>
&=&
\<\phi,\mu\> + t\<\phi,m\>
+ \frac{1}{2}\int_0^t\<a\phi^{\prime\prime}, X_s\> ds
- \int_0^t\<b\phi,X_s\> ds \nonumber \\
& & + \int_0^t\int_{\mbb{R}}\phi(y) Z(ds,dy) +
\int_0^t\int_{\mbb{R}}\<h(y-\cdot) \phi^\prime,X_s\> W(ds,dy),
 \end{eqnarray}
where $W(ds,dx)$ is a time-space white noise and $Z(ds,dy)$ is an
orthogonal martingale measure which is orthogonal to $W(ds,dy)$
and has covariation measure $\sigma(y)X_s(dy)ds$. Let $({\cal
F}_t)_{t\ge0}$ denote the filtration generated by $\{W(ds,dy)\}$
and $\{Z(ds,dy)\}$. Since $\sigma$ is strictly positive, the
process $\{X_t: t\ge0\}$ can be represented in terms of the
covariation measure of $Z(ds,dy)$, so it is adapted to $({\cal
F}_t)_{t\ge0}$. By Theorem~\ref{t4.2}, for $\phi\in H_1(\mbb{R})
\cap C_b(\mbb{R})^+$ the equation
 \begin{eqnarray}\label{5.2}
\psi_{r,t}(x)
&=&
\phi(x) + \int_r^t \bigg[\frac{1}{2} a(x)\psi_{s,t}
^{\prime\prime}(x) - b(x)\psi_{s,t}(x) - \frac{1}{2}\sigma(x)
\psi_{s,t}(x)^2\bigg] ds \nonumber \\
& & + \int_r^t\int_{\mbb{R}} h(y-x)\psi_{s,t}^\prime(x) \cdot
W(ds,dy), \qquad t\ge r\ge 0,
 \end{eqnarray}
has a pathwise unique solution $\psi_{r,t} = \psi^W_{r,t}$ in
$H_1(\mbb{R}) \cap C_b(\mbb{R})^+$. Let $\mbf{P}^W$ and
$\mbf{E}^W$ denote respectively the conditional probability and
expectation given the white noise $\{W(ds,dy)\}$. The main result
of this section is the following

\btheorem\label{t5.1} For $t\ge r\ge0$ and $\phi\in H_1(\mbb{R})
\cap C_b(\mbb{R})^+$ we have a.s.\
 \begin{eqnarray}\label{5.3}
\mbf{E}^W\{ e^{-\<\phi,X_t\>}|{\cal F}_r\} =
\exp\bigg\{-\<\psi^W_{r,t},X_r\> -
\int_r^t\<\psi^W_{s,t},m\>ds\bigg\},
 \end{eqnarray}
where $\psi^W_{r,t}$ is defined by (\ref{5.2}). Consequently,
$\{X_t: t\ge0\}$ is a diffusion process with Feller transition
semigroup $(Q_t)_{t\ge 0}$ given by
 \begin{eqnarray}\label{5.4}
\int_{M(\mbb{R})} e^{-\<\phi,\nu\>} Q_t(\mu,d\nu) =
\mbf{E}\exp\bigg\{-\<\psi^W_{0,t},\mu\> -
\int_0^t\<\psi^W_{s,t},m\>ds\bigg\}.
 \end{eqnarray}
\etheorem

Our proof of the theorem are based on direct calculations derived
from (\ref{5.1}) and (\ref{5.2}). The argument is different from
that of \cite{X04}, where the Wong-Zakai approximation was used to
get the result. We shall give four lemmas which together with the
proof of the theorem show clearly the key steps of the
calculations.

Suppose that $\alpha$ and $\beta$ are bounded measurable functions
on $[0,\infty) \times \mbb{R}$ and that
 \begin{eqnarray*}
\int_0^t\int_{\mbb{R}} \alpha(s,y)^2dsdy < \infty.
 \end{eqnarray*}
For $t\ge r\ge0$, define
 \begin{eqnarray}\label{5.5}
\theta_\alpha(r,t) = \exp\bigg\{\int_r^t\int_{\mbb{R}}
\alpha(s,y)W(ds,dy) - \frac{1}{2}\int_r^t\int_{\mbb{R}}
\alpha(s,y)^2dsdy\bigg\},
 \end{eqnarray}
and
 \begin{eqnarray}\label{5.6}
\zeta_\beta(r,t) = \exp\bigg\{\int_r^t\int_{\mbb{R}}
\beta(s,y)Z(ds,dy) - \frac{1}{2}\int_r^t \<\sigma
\beta(s,\cdot)^2,X_s\> ds\bigg\}.
 \end{eqnarray}
Then we have the following

\blemma\label{l5.1} Under the conditional probability measure
$\mbf{P}^W$, the process $\{\zeta_\beta(0,t): t\ge 0\}$ is a
martingale with respect to $({\cal F}_t)_{t\ge 0}$. \elemma

\proof Clearly, both $\{\theta_\alpha(0,t): t\ge 0\}$ and
$\{\zeta_\beta(0,t): t\ge 0\}$ are martingales under the original
probability measure $\mbf{P}$. Recall that the martingale measures
$\{W(ds,dy)\}$ and $\{Z(ds,dy)\}$ are orthogonal. By integration
by parts it is easy to see that $\{\theta_\alpha(0,t)
\zeta_\beta(0,t): t\ge 0\}$ is a martingale. Since $\alpha$ is
arbitrary, for any $u\ge t\ge r\ge0$ and any bounded ${\cal
F}_r$-measurable random variable $Z$ we obtain
 \begin{eqnarray*}
\mbf{E}\{\theta_\alpha(0,u)\zeta_\beta(0,t)Z\} =
\mbf{E}\{\theta_\alpha(0,r)\zeta_\beta(0,r)Z\} =
\mbf{E}\{\theta_\alpha(0,u)\zeta_\beta(0,r)Z\}.
 \end{eqnarray*}
Note that the linear span of the functionals $\{\theta_\alpha
(0,u)\}$ is dense in the space of squared-integrable and
$\sigma(W)$-measurable random variables; see e.g.\
\cite[p.81]{B92} and \cite{C04}. Then we have the desired equality
$\mbf{E}^W \{\zeta_\beta(0,t)| {\cal F}_r\} = \zeta_\beta(0,r)$.
\qed

By the property of independent increments of the white noise
$\{W(ds,dy)\}$ we have
 \begin{eqnarray}\label{5.7}
\xi_{r,t}(x) := \mbf{E}\{\psi_{r,t}(x)\theta_\alpha(r,t)\} =
\mbf{E}\{\psi_{r,t}(x)\theta_\alpha(r,t)|{\cal F}_r\}
 \end{eqnarray}
and
 \begin{eqnarray}\label{5.8}
\eta_{r,t}(x) := \mbf{E}\{\psi_{r,t}(x)^2\theta_\alpha(r,t)\} =
\mbf{E}\{\psi_{r,t}(x)^2\theta_\alpha(r,t)|{\cal F}_r\}.
 \end{eqnarray}

\blemma\label{l5.2} For $t\ge r\ge 0$, we have a.s.\
 \begin{eqnarray}\label{5.9}
\mbf{E}\{\<\psi_{r,t},X_r\>\theta_\alpha(0,t)\zeta_\beta(0,t)|{\cal
F}_r\} = \<\xi_{r,t},X_r\>\theta_\alpha(0,r)\zeta_\beta(0,r)
 \end{eqnarray}
and
 \begin{eqnarray}\label{5.10}
\mbf{E}\{\<\sigma\psi_{r,t}^2,X_r\>\theta_\alpha(0,t)\zeta_\beta(0,t)
| {\cal F}_r\} =
\<\sigma\eta_{r,t},X_r\>\theta_\alpha(0,r)\zeta_\beta(0,r).
 \end{eqnarray}
\elemma

\proof By Lemma~\ref{l5.1} it is easy to see that
$\mbf{E}^W[\zeta_\beta(r,t) |{\cal F}_r] = 1$. Since
$\theta_\alpha (0,r) \zeta_\beta(0,r)$ is ${\cal F}_r$-measurable
and $\<\psi_{r,t},X_r\> \theta_\alpha(r,t)$ is $\sigma(W, {\cal
F}_r)$-measurable, we have
 \begin{eqnarray*}
&& \mbf{E}\{\<\psi_{r,t},X_r\>\theta_\alpha(0,t)
\zeta_\beta(0,t)|{\cal F}_r\} \nonumber \\
 &&\qquad=
\mbf{E}\{\<\psi_{r,t},X_r\>\theta_\alpha(r,t)\zeta_\beta(r,t)|{\cal
F}_r\}
\theta_\alpha(0,r)\zeta_\beta(0,r) \nonumber \\
 &&\qquad=
\mbf{E}\{\<\psi_{r,t},X_r\>\theta_\alpha(r,t)\mbf{E}^W[\zeta_\beta(r,t)
|{\cal F}_r]|{\cal F}_r\}\theta_\alpha(0,r)\zeta_\beta(0,r) \nonumber \\
 &&\qquad=
\mbf{E}\{\<\psi_{r,t},X_r\>\theta_\alpha(r,t)
|{\cal F}_r\}\theta_\alpha(0,r)\zeta_\beta(0,r) \nonumber \\
 &&\qquad=
\<\xi_{r,t},X_r\>\theta_\alpha(0,r)\zeta_\beta(0,r).
 \end{eqnarray*}
A similar calculation gives (\ref{5.10}). \qed

\blemma\label{l5.3} For $t\ge r\ge 0$ and $x\in \mbb{R}$, we have
 \begin{eqnarray}\label{5.11}
\xi_{r,t}(x)
 &=&
\phi(x) + \int_r^t
\bigg[\frac{1}{2}a(x)\xi_{s,t}^{\prime\prime}(x) - b(x)
\xi_{s,t}(x) - \frac{1}{2}\sigma(x)\eta_{s,t}(x)\bigg]
ds \nonumber \\
 & &
+ \int_r^t\<h(\cdot-x), \alpha(s,\cdot)\> \xi_{s,t}^\prime(x)ds,
 \end{eqnarray}
where the derivatives are taken in the classical sense. \elemma

\proof Note that the backward and forward integrals coincide for
deterministic integrands. Then we may fix $t>0$ and apply It\^o's
formula to the process $\{\theta_\alpha(r,t): r\in [0,t]\}$ to get
 \begin{eqnarray}\label{5.12}
\theta_\alpha(r,t) = 1 + \int_r^t\int_{\mbb{R}} \theta_\alpha(s,t)
\alpha(s,y) \cdot W(ds,dy).
 \end{eqnarray}
By (\ref{5.2}), (\ref{5.12}) and backward It\^o formula, for any
$f\in C^\infty_b(\mbb{R})$ we have
 \begin{eqnarray}\label{5.13}
\<\psi_{r,t},f\>\theta_\alpha(r,t)
 &=&
\<\phi,f\> + \int_r^t\bigg[\frac{1}{2}\<a\psi_{s,t}
^{\prime\prime}, f\> - \<b\psi_{s,t},f\> - \frac{1}{2}
\<\sigma\psi_{s,t}^2,f\>\bigg]\theta_\alpha(s,t) ds \nonumber \\
 & &
+ \int_r^t\int_{\mbb{R}} [\<h(y-\cdot)\psi_{s,t}^\prime,f\> +
\<\psi_{s,t},f\>\alpha(s,y)]\theta_\alpha(s,t) \cdot
W(ds,dy) \nonumber \\
 & &
+ \int_r^t\int_{\mbb{R}} \<h(y-\cdot)\psi_{s,t}^\prime,f\>
\theta_\alpha(s,t)\alpha(s,y) ds dy.
 \end{eqnarray}
(See e.g.\ \cite[p.124]{B92} for the backward It\^o formula.)
Observe that for fixed $t>0$, the process
 \begin{eqnarray*}
\int_r^t\int_{\mbb{R}} [\<h(y-\cdot)\psi_{s,t}^\prime,f\> +
\<\psi_{s,t},f\>\alpha(s,y)]\theta_\alpha(s,t) \cdot W(ds,dy)
 \end{eqnarray*}
is a backward martingale in $r\le t$. Taking the
expectation in (\ref{5.13}) we obtain
 \begin{eqnarray*}
\<\xi_{r,t},f\>
 &=&
\<\phi,f\> + \int_r^t \bigg[\frac{1}{2}
\<a\xi_{s,t}^{\prime\prime}, f\> - \<b\xi_{s,t},f\>
- \frac{1}{2} \<\sigma\eta_{s,t},f\>\bigg] ds \nonumber \\
& & + \int_r^t\int_{\mbb{R}} \<h(y-\cdot)\xi_{s,t}^\prime,f\>
a(s,y) ds dy.
 \end{eqnarray*}
Then $\{\xi_{r,t}\}$ must coincides with the classical solution of
the parabolic equation (\ref{5.11}). \qed

\blemma\label{l5.4} For any $t\ge r\ge 0$, we have a.s.\
 \begin{eqnarray}\label{5.14}
\<\phi,X_t\> = \<\psi_{r,t},X_r\> +
\int_r^t\int_{\mbb{R}}\psi_{s,t}(x)Z(ds,dx) + \frac{1}{2}\int_r^t
\<\sigma\psi_{s,t}^2,X_s\>ds + \int_r^t \<\psi_{s,t},m\>ds.
 \end{eqnarray}
\elemma

\proof In view of (\ref{5.1}) and (\ref{5.11}), we may integrate
$\xi_{s,t}$ backward relative to $X_s$ to see that
 \begin{eqnarray*}
d\<\xi_{s,t},X_s\>
 &=&
\frac{1}{2}\<\sigma\eta_{s,t},X_s\> ds - \int_{\mbb{R}}
\<h(y-\cdot)\xi_{s,t}^\prime,X_s\> \alpha(s,y) ds dy
+ \<\xi_{s,t},m\>ds \nonumber \\
 & &
+ \int_{\mbb{R}} \xi_{s,t}(y) Z(ds,dy) + \int_{\mbb{R}}
\<h(y-\cdot) \xi_{s,t}^\prime,X_s\> W(ds,dy),
 \end{eqnarray*}
where the first two terms from (\ref{5.11}) cancelled out with the
second and third terms from (\ref{5.1}). Since the two martingale
measures $\{W(ds,dy)\}$ and $\{Z(ds,dy)\}$ are orthogonal, by
It\^o's formula we have
 \begin{eqnarray}\label{5.15}
d\<\xi_{s,t},X_s\> \theta_\alpha(0,s)\zeta_\beta(0,s)
 &=&
\frac{1}{2}\<\sigma\eta_{s,t},X_s\>\theta_\alpha(0,s)
\zeta_\beta(0,s) ds + \<\xi_{s,t},m\>\theta_\alpha(0,s)
\zeta_\beta(0,s)ds  \nonumber \\
 & &
+ \int_{\mbb{R}} \<h(y-\cdot)\xi_{s,t}^\prime,X_s\>
\theta_\alpha(0,s)\zeta_\beta(0,s) W(ds,dy)  \nonumber \\
 & &
+ \int_{\mbb{R}} \xi_{s,t}(y)\theta_\alpha(0,s)
\zeta_\beta(0,s)Z(ds,dy)   \nonumber \\
 & &
+ \int_{\mbb{R}} \<\xi_{s,t},X_s\>\theta_\alpha(0,s)
\zeta_\beta(0,s)\alpha(s,x) W(ds,dy)  \nonumber \\
 & &
+ \int_{\mbb{R}} \<\xi_{s,t},X_s\>\theta_\alpha(0,s)
\zeta_\beta(0,s)\beta(s,y) Z(ds,dy)  \nonumber \\
 & &
+ \,\<\sigma\xi_{s,t} \beta(s,\cdot),X_s\>
\theta_\alpha(0,s)\zeta_\beta(0,s)ds.
 \end{eqnarray}
By a calculation similar to the proof of Lemma~\ref{l5.2} we get
 \begin{eqnarray}\label{5.16}
\mbf{E}\{\<\psi_{s,t},m\>\theta_\alpha(0,t)\zeta_\beta(0,t)|{\cal
F}_s\} = \<\xi_{s,t},m\>\theta_\alpha(0,s)\zeta_\beta(0,s).
 \end{eqnarray}
{From} (\ref{5.10}), (\ref{5.15}) and (\ref{5.16}) it follows that
 \begin{eqnarray}\label{5.17}
 & &\mbf{E}\{\<\phi,X_t\>\theta_\alpha(0,t)\zeta_\beta(0,t)\}
- \mbf{E}\{\<\xi_{r,t},X_r\>\theta_\alpha(0,r)
\zeta_\beta(0,r)\} \nonumber \\
 &=&
\frac{1}{2}\mbf{E}\bigg\{\int_r^t\<\sigma\eta_{s,t},X_s\>
\theta_\alpha(0,s)\zeta_\beta(0,s) ds\bigg\} +
\mbf{E}\bigg\{\int_r^t\<\xi_{s,t},m\>
\theta_\alpha(0,s)\zeta_\beta(0,s)ds\bigg\} \nonumber \\
 & &
+\, \mbf{E}\bigg\{\int_r^t\<\sigma\xi_{s,t} \beta(s,\cdot),X_s\>
\theta_\alpha(0,s)\zeta_\beta(0,s)ds\bigg\}  \nonumber \\
 &=&
\frac{1}{2}\mbf{E}\bigg\{\int_r^t\<\sigma\psi_{s,t}^2,X_s\>
\theta_\alpha(0,t)\zeta_\beta(0,t) ds\bigg\} +
\mbf{E}\bigg\{\int_r^t\<\psi_{s,t},m\> \theta_\alpha(0,t)
\zeta_\beta(0,t)ds\bigg\} \nonumber \\
 & &
+\, \mbf{E}\bigg\{\int_r^t\<\sigma\xi_{s,t} \beta(s,\cdot),X_s\>
\theta_\alpha(0,s)\zeta_\beta(0,s)ds\bigg\}.
 \end{eqnarray}
By (\ref{5.6}) and It\^o's formula we have
 \begin{eqnarray*}
\zeta_\beta(0,t) = 1 + \int_0^t\int_{\mbb{R}}
\zeta_\beta(0,s)\beta(s,y) Z(ds,dy),
 \end{eqnarray*}
and hence
 \begin{eqnarray*}
 &&\mbf{E}\bigg\{\int_r^t\int_{\mbb{R}}\psi_{s,t}(y) Z(ds,dy)
\theta_\alpha(0,t)\zeta_\beta(0,t)\bigg\}   \\
 &&\qquad =
\mbf{E}\bigg\{\mbf{E}^W\bigg[\int_r^t\int_{\mbb{R}} \psi_{s,t}(y)
Z(ds,dy)\zeta_\beta(0,t)\bigg]
\theta_\alpha(0,t)\bigg\}   \\
 &&\qquad =
\mbf{E}\bigg\{\mbf{E}^W\bigg[\int_r^t\<\sigma
\psi_{s,t}\beta(s,\cdot),X_s\> \zeta_\beta(0,s)ds\bigg]
\theta_\alpha(0,t)\bigg\}   \\
 &&\qquad =
\int_r^t\mbf{E}\bigg[\<\sigma
\psi_{s,t}\beta(s,\cdot),X_s\>\theta_\alpha(0,t)
\zeta_\beta(0,s)\bigg]ds   \\
 &&\qquad =
\int_r^t\mbf{E}\bigg[\<\sigma
\xi_{s,t}\beta(s,\cdot),X_s\>\theta_\alpha(0,s)
\zeta_\beta(0,s)\bigg]ds,
 \end{eqnarray*}
where the last equality follows from (\ref{5.9}). Then we
substitute the above into (\ref{5.17}) to get
 \begin{eqnarray*}
 && \mbf{E}\{\<\phi,X_t\>\theta_\alpha(0,t)\zeta_\beta(0,t)\} -
\mbf{E}\{\<\xi_{r,t},X_r\>\theta_\alpha(0,r)
\zeta_\beta(0,r)\} \nonumber \\
 &&\qquad =
\frac{1}{2}\mbf{E}\bigg\{\int_r^t\<\sigma\psi_{s,t}^2,X_s\>
\theta_\alpha(0,t)\zeta_\beta(0,t) ds\bigg\} +
\mbf{E}\bigg\{\int_r^t\<\psi_{s,t},m\> \theta_\alpha(0,t)
\zeta_\beta(0,t)ds\bigg\} \nonumber \\
 & &\qquad\qquad
+\, \mbf{E}\bigg\{\int_r^t\int_{\mbb{R}}\psi_{s,t}(y)
Z(ds,dy)\theta_\alpha(0,t)\zeta_\beta(0,t)\bigg\}.
 \end{eqnarray*}
On the other hand, by (\ref{5.9}) we have
 \begin{eqnarray*}
& &\mbf{E}\{\<\phi,X_t\>\theta_\alpha(0,t)\zeta_\beta(0,t)\}
- \mbf{E}\{\<\xi_{r,t},X_r\>\theta_\alpha(0,r)\zeta_\beta(0,r)\}   \\
 &=&
\mbf{E}\{[\<\phi,X_t\> - \<\psi_{r,t},X_r\>]
\theta_\alpha(0,t)\zeta_\beta(0,t)\}.
 \end{eqnarray*}
It follows that
 \begin{eqnarray*}
\lefteqn{\mbf{E}\bigg\{\bigg[\<\phi,X_t\> - \<\psi_{r,t},X_r\> -
\frac{1}{2}\int_r^t \<\sigma\psi_{s,t}^2,X_s\> ds}  \\
 & &
- \int_r^t \<\psi_{s,t},m\>ds -
\int_r^t\int_{\mbb{R}}\psi_{s,t}(x)Z(ds,dx)\bigg]
\theta_\alpha(0,t)\zeta_\beta(0,t)\bigg\} = 0.
 \end{eqnarray*}
Then we have the desired equation; see e.g.\ \cite[p.81]{B92} and
\cite{C04}. \qed

\noindent\textit{Proof of Theorem~\ref{t5.1}.~} Recall that
$Z(ds,dy)$ is an orthogonal martingale measure with covariation
measure $\sigma(y) X_s(dy)ds$. By Lemma~\ref{l5.1}, for any fixed
$u\ge r$ the process
 \begin{eqnarray*}
\exp\bigg\{- \int_r^t\int_{\mbb{R}}\psi_{s,u}(y)Z(ds,dy) -
\frac{1}{2}\int_r^t \<\sigma\psi_{s,u}^2,X_s\> ds\bigg\}, \qquad
r\le t\le u,
 \end{eqnarray*}
is a martingale under $\mbf{P}^W$. By Lemma~\ref{l5.4} we get
a.s.\
 \begin{eqnarray*}
\mbf{E}^W\{ e^{-\<\phi,X_t\>}|{\cal F}_r\}
 &=&
\mbf{E}^W\bigg[ \exp\bigg\{-\<\psi_{r,t},X_r\>
- \int_r^t\int_{\mbb{R}}\psi_{s,t}(y)Z(ds,dy)    \\
 & & - \frac{1}{2}\int_r^t \<\sigma\psi_{s,t}^2,
X_s\> ds - \int_r^t\<\psi_{s,t},m\>ds\bigg\}
\bigg|{\cal F}_r\bigg] \nonumber  \\
&=& \exp\bigg\{-\<\psi_{r,t},X_r\>
- \int_r^t\<\psi_{s,t},m\>ds\bigg\},
 \end{eqnarray*}
giving (\ref{5.3}). In particular, we have
 \begin{eqnarray}\label{5.18}
\mbf{E}\{ e^{-\<\phi,X_t\>}\} =
\mbf{E}\exp\bigg\{-\<\psi_{0,t},\mu\> -
\int_0^t\<\psi_{s,t},m\>ds\bigg\}.
 \end{eqnarray}
The distribution of $X_t$ is uniquely determined by (\ref{5.18})
and the uniqueness of solution of (\ref{5.1}) follows. This in
turn implies the strong Markov property of $\{X_t: t\ge0\}$. Since
$\psi_{r,t} (x)$ is continuous in $x\in \mbb{R}$, the transition
semigroup $(Q_t)_{t\ge 0}$ defined by (\ref{5.4}) is Feller. \qed


\section{Some properties of the SDSMI}

\setcounter{equation}{0}

We here investigate some properties of the SDSMI. Let $(c, h,
\sigma, b, m)$ be given as in the introduction. As in the last
section, let $\mbf{P}^W$ and $\mbf{E}^W$ denote respectively the
conditional probability and expectation given the white noise
$\{W(ds,dy)\}$. The equality (\ref{5.3}) suggests that $\{X_t:
t\ge 0\}$ under $\mbf{P}^W$ is a Markov process with transition
semigroup $(Q^W_{r,t})_{t\ge r}$ satisfying a.s.\
 \begin{eqnarray}\label{6.1}
\int_{M(E)}e^{-\<\phi,\nu\>} Q^W_{r,t}(\mu, d\nu) =
\exp\bigg\{-\<\psi^W_{r,t},\mu\> -
\int_r^t\<\psi^W_{s,t},m\>ds\bigg\}.
 \end{eqnarray}
In other words, the SDSMI conditioned upon $\{W(ds,dy)\}$ should
be an inhomogeneous immigration superprocess. This observation
suggests a number of applications of the conditional log-Laplace
functional. For instance, based on the results in the last
section, the conditional excursion theory of the SDSM have been
developed in \cite{LWX04b}. Moreover, some moment formulas can be
also derived from (\ref{5.3}) in a similar way as \cite{X04}.

As another application of the conditional Laplace functionals, we
prove the following ergodicity property of the SDSMI.

\btheorem\label{t6.1} Suppose that there is a constant $\epsilon
>0$ such that $b(x) \ge \epsilon$ for all $x\in \mbb{R}$. Then the
SDSMI has a unique stationary distribution $Q_\infty$ given by
 \begin{eqnarray}\label{6.2}
\int_{M(\mbb{R})} e^{-\<\phi,\nu\>} Q_\infty(d\nu) =
\mbf{E}\exp\bigg\{- \int_0^\infty\<\psi^W_t,m\>dt\bigg\},
 \end{eqnarray}
where $\psi^W_t(x)$ is the solution of (\ref{4.1}). Moreover, we
have $\lim_{t\to\infty} Q_t(\mu,\cdot) = Q_\infty(\cdot)$ in the
topology of weak convergence for each $\mu \in M(\mbb{R})$.
\etheorem

\proof Using the notation of the proof of Theorem~\ref{t4.2}, for
any $t\ge r\ge 0$ we have
 \begin{eqnarray*}
\mbf{E}\exp\bigg\{-\int_r^t\<\psi^W_{s,t},m\>ds\bigg\}
 &=&
\mbf{E}\exp\bigg\{-\int_r^t\<\phi^W_{t-s,t},m\>ds\bigg\}  \\
 &=&
\mbf{E}\exp\bigg\{-\int_0^{t-r}\<\phi^W_{s,t},m\>ds\bigg\}  \\
 &=&
\mbf{E}\exp\bigg\{-\int_0^{t-r}\<\psi^W_s,m\>ds\bigg\},
 \end{eqnarray*}
where the last equality follows by the property of independent and
stationary increments of the time-space white noise. By
Theorem~\ref{t4.2} we have $\|\psi^W_{s,t}\| \le e^{-\epsilon
(t-s)}\|\phi\|$ for $s\le t$. It follows that
 \begin{eqnarray*}
\lim_{t\to\infty}\int_{M(\mbb{R})} e^{-\<\phi,\nu\>} Q_t(\mu,d\nu)
&=& \lim_{t\to\infty}\mbf{E}\exp\bigg\{-\<\psi^W_{0,t},\mu\>
- \int_0^t\<\psi^W_{s,t},m\>ds\bigg\}  \\
&=& \lim_{t\to\infty}\mbf{E}\exp\bigg\{
- \int_0^t\<\psi^W_{s,t},m\>ds\bigg\}  \\
&=& \mbf{E}\exp\bigg\{ - \int_0^{\infty}\<\psi^W_s,m\>ds\bigg\}.
 \end{eqnarray*}
On the other hand, by Theorem~\ref{t4.1} it is easy to get
 \begin{eqnarray*}
\lim_{\|\phi\|\to 0}\mbf{E}\exp\bigg\{- \int_0^{\infty}
\<\psi^W_s,m\>ds\bigg\} = 1.
 \end{eqnarray*}
Then (\ref{6.2}) defines a probability measure $Q_\infty$ on
$M(\mbb{R})$ and $\lim_{t\to\infty} Q_t(\mu,\cdot) =
Q_\infty(\cdot)$ in the topology of weak convergence; see e.g.\
\cite[Lemma~2.1]{L02}. \qed

The properties of the SDSMI varies sharply for different choices
of the parameters. The special case where $b(\cdot) \equiv 0$ and
$\<1, m\> = 0$ was discussed in \cite{DLW01, DVW00, W97, W02}. In
this case, we have
 \begin{eqnarray}\label{6.3}
\<\phi,X_t\> &=& \<\phi,\mu\>
+ \frac{1}{2}\int_0^t\<a\phi^{\prime\prime}, X_s\> ds
+ \int_0^t\int_{\mbb{R}}\phi(y) Z(ds,dy) \nonumber \\
& & + \int_0^t\int_{\mbb{R}}\<h(y-\cdot) \phi^\prime,X_s\>
W(ds,dy).
 \end{eqnarray}
The solution of (\ref{6.3}) is a critical branching SDSM without
immigration. In particular, if $c(\cdot)$ is bounded away from
zero, then $\{X_t: t>0\}$ is absolutely continuous for any initial
state $X_0$; see \cite{DLW01, DVW00, W97}. On the other hand, if
$c(\cdot) \equiv 0$, then $\{X_t: t>0\}$ is purely atomic for any
initial state $X_0$; see \cite{DL03, W97, W02}.

Another special case is where $\sigma(\cdot) \equiv 0$ and
$\<1,m\> = 0$. In this case, we get from (\ref{6.3}) the linear
equation
 \begin{eqnarray}\label{6.4}
\<\phi,X_t\> = \<\phi,\mu\> +
\frac{1}{2}\int_0^t\<a\phi^{\prime\prime}, X_s\> ds -
\int_0^t\<b\phi,X_s\> ds + \int_0^t\int_{\mbb{R}}\<h(y-\cdot)
\phi^\prime,X_s\> W(ds,dy).
 \end{eqnarray}
The process defined in this way is closely related to the
superprocesses arising from isotropic stochastic flows
investigated by \cite{MX01}. The following theorem shows that
$\{X_t: t\ge0\}$ is absolutely continuous for a large class of
absolutely continuous initial states.

\btheorem\label{t6.2} If $\{X_t: t\ge0\}$ is a solution of
(\ref{6.4}) with $X_0(dx) = v_0(x)dx$ for some $v_0\in
H_0(\mbb{R})$, then there is an $H_0(\mbb{R})$-valued process
$\{v_t: t\ge0\}$ such that $X_t(dx) = v_t(x)dx$ a.s.\ holds.
\etheorem

\proof By \cite[Theorem~3.5]{KX99}, the equation
 \begin{eqnarray}\label{6.5}
v_t(x) = v_0(x) + \int_0^t
\bigg[\frac{1}{2}(av_s)^{\prime\prime}(x) - b(x)v_s(x)\bigg] ds -
\int_0^t\int_{\mbb{R}} (h(y-\cdot)v_s)^\prime(x) W(ds,dy)
 \end{eqnarray}
has a unique $H_0(\mbb{R})$-valued solution $\{v_t: t\ge0\}$. Let
$X_t(dx) = v_t(x)dx$. Clearly, $\{X_t: t\ge0\}$ solves
(\ref{6.4}). \qed

\end{document}